%% file: extremal2a.tex
\documentclass{article}

\begin{document}

\title{Constant scalar curvature  metrics on toric surfaces}
\author{S. K. Donaldson}

\maketitle

\tableofcontents

\newtheorem{thm}{Theorem}
\newtheorem{prop}{Proposition}
\newtheorem{defn}{Definition}
\newtheorem{condn}{Condition}
\newtheorem{cor}{Corollary}
\newtheorem{lem}{Lemma}
\newcommand{\tu}{\tilde{u}}
\newcommand{\up}{\underline{p}}
\newcommand{\uy}{\underline{y}}
\newcommand{\ux}{{\underline{x}}}
\newcommand{\utu}{\underline{\tilde{u}}}
\newcommand{\uh}{\underline{h}}
\newcommand{\oh}{\overline{h}}
\newcommand{\oP}{\overline{P}}
\newcommand{\tP}{\tilde{P}}
\newcommand{\tA}{\tilde{A}}
\newcommand{\ut}{\underline{t}}
\newcommand{\ot}{\overline{t}}
\newcommand{\tsigma}{\tilde{\sigma}}
\newcommand{\tU}{\tilde{U}}
\newcommand{\ufl}{u_{\flat}}
\newcommand{\xifl}{\xi^{\flat}}
\newcommand{\ua}{\underline{a}}
\newcommand{\Pfl}{P_{\flat}}
\newcommand{\Afl}{A_{\flat}}
\newcommand{\Vfl}{V_{\flat}}
\newcommand{\Area}{{\rm Area}\ }
\newcommand{\oM}{\overline{M}}
\newcommand{\bR}{{\bf R}}
\newcommand{\bC}{{\bf C}}
\newcommand{\uu}{\underline{u}}
\newcommand{\Vol}{{\rm Vol}}
\newcommand{\Emax}{E_{{\rm max}}}
\newcommand{\Fmax}{F_{{\rm max}}}
\newcommand{\udelta}{\underline{\delta}}
\newcommand{\xist}{\xi^{*}}
\newcommand{\distu}{{\rm dist}_{u}}
\newcommand{\distual}{{\rm dist}_{u^{(\alpha)}}}

\section{Introduction}

This paper continues the series \cite{kn:D1}, \cite{kn:D2}, \cite{kn:D3}
in which we study the scalar curvature of Kahler metrics on toric varieties and relations with the analysis of convex
functions on  polytopes in Euclidean space. The main result
of the present paper is an existence theorem for metrics of constant scalar
curvature on toric surfaces, confirming a conjecture in \cite{kn:D1}. 
 
       We begin by recalling the background briefly: more details can be
       found in the references above.  Let $P$ be a bounded open polytope in $\bR^{n}$ and let $\sigma$      be
        a measure on the boundary of $P$ which is a multiple of the standard
        Lebesgue measure on each face. Let $A$ be a smooth function on the
        closure $\oP$ of $P$ and consider the linear functional $L_{A,\sigma}$
        on the continuous functions on $\oP$ given by
        $$  L_{A,\sigma} f= \int_{\partial P} f\  d\sigma - \int_{P} A f \ d\mu,
        $$
        where $d\mu$ is ordinary Lebesgue measure on $\bR^{n}$. We define
        a nonlinear functional  on a suitable class of {\it convex}  functions $u$ on $P$ by
        
           $${\cal M}(u)= - \int_{P} \log \det (u_{ij})\ d\mu + L_{A,\sigma}
           u, $$
           where $(u_{ij})$ denotes the Hessian matrix of second derivatives of $u$.
           
             The Euler-Lagrange equation associated to ${\cal M}$ is the
             fourth order PDE found by Abreu:
             \begin{equation}\sum_{i, j} \frac{\partial^{2} u^{ij}}{\partial x_{i} \partial
             x_{j}}  = -A, \end{equation}
             where $(u^{ij})$ is the matrix inverse of $(u_{ij})$. (We often
             use the notation $u^{ij}_{ij}$ for the left-hand side of (1).) A solution
             of this equation, with appropriate boundary  behaviour, is a critical
             point, in fact the minimiser, of ${\cal M}$. More precisely,
             we require $u$ to satisfy  {\it Guillemin boundary conditions},
             which depend on the measure $\sigma$. We refer to the previous
             papers cited for the details of these boundary conditions and
             just recall here the  standard model for the behaviour
             at the boundary. This
             is the function 
             $$u(x_{1}, \dots, x_{n}) = \sum x_{i} \log x_{i}, $$
             on the convex subset $\{ x_{i} >0\}$ in $\bR^{n}$. For appropriate
             \lq\lq Delzant'' pairs  $(P,\sigma)$ a function $u$ with this
             boundary behaviour defines a Kahler metric on a corresponding
             toric manifold $X_{P}$. This  comes with a map $\pi:X_{P}\rightarrow
             P$ and the scalar curvature of the metric is $A\circ \pi$. Thus
             when $A$ is constant the Kahler metric has constant scalar curvature.
             When $A$ is an affine-linear function the Kahler metric is \lq\lq
             extremal''.
             
              The Guillemin
              boundary conditions imply that
             \begin{equation} L_{A,\sigma} (f)= \int_{P} \sum f_{ij} u^{ij} d\mu, \end{equation} for
            smooth  test functions $f$. Thus $L_{A}$ vanishes on affine-linear functions.
             This  just says that  $(P,A d\mu)$ and $(\partial P, d\sigma)$ have the same mass and moments. More interestingly, equation (2) tells us
that $L_{A,\sigma}(f)\geq 0$ for  smooth {\it convex} functions $f$, with strict
equality if $f$ is not affine-linear. 
            This can be
    extended to more general convex functions $f$.
    In \cite{kn:D1} we conjectured that this necessary
    condition for the existence of a solution $u$, is also sufficient. The
    present paper completes the proof of this in the case when the dimension
    $n$ is $2$ and the function $A$ is a constant. Thus we have
    \begin{thm}
    Suppose $P\subset \bR^{2}$ is a polygon and  $\sigma$ is a measure on $\partial
    P$, as above, with the property that the mass and moments of $(\partial
    P, \sigma)$ and $(P, A d\mu)$ are equal for some constant $A$. Then {\it
    either} there is a solution to (1), satisfying Guillemin boundary conditions,
   {\it  or} there is convex function $f$, not affine-linear, with
   $L_{A,\sigma}(f)\leq 0$
   \end{thm}

    A corollary  in the framework of complex geometry of this is
    
    \begin{cor}
    If a polarised complex toric surface with zero Futaki invariant is K-stable
    it admits a constant
    scalar curvature Kahler metric. 
    \end{cor}
    We refer to \cite{kn:D1} for the terminology and further details. (The
    converse is also true and is proved by Zhou and Zhu in \cite{kn:ZZ}.)

   \
   
   \ Two remarks are in order here.
    First, in the
    two-dimensional case the 
     \lq\lq positivity'' condition for the functional $L_{A}$ can be made
     more explicit. If $\lambda$ is an affine-linear function on $\bR^{2}$
     define
     $$  \lambda^{+}(\ux)= \max(0,\lambda(\ux)). $$
     Then it follows from the arguments of \cite{kn:D1} that the second alternative
     in Theorem is equivalent to the condition that there is some function
     of the form $\lambda^{+}$, not identically zero on $P$, with $L_{A,\sigma}(\lambda^{+})\leq
     0$. Thus to check the existence of a solution to (1) we only need to
     check the positivity of $L_{A,\sigma}$ on a $2$-parameter family of         functions
     of the form $\lambda^{+}$, which would be easy to do with a computer.
     Second, it was shown in \cite{kn:D1} that the positivity condition is
     equivalent to a more useful, quantitative, statement. 
      Fix a base point 
    $p_{0}$ and say a convex function $f$ is {\it normalised at $p_{0}$} if
    $p_{0}$ is the minimiser of $f$ and $f(p_{0})=0$. Then the positivity property holds if and only if there is some $C= C_{A,\sigma, p_{0}}>0$ such that
   
    \begin{equation}  \int_{\partial P} f d\sigma \leq C \ L_{A,\sigma} f \end{equation}
    for all normalised convex functions $f$. 
    \

    The main result of \cite{kn:D3} asserts that Theorem 1 follows if one
    can establish a certain {\it a priori} estimate on solutions $u$. Suppose that
    (3) holds and $u$ is a solution, normalised at $p_{0}$ Then taking $f=u$ in (2) we obtain
    \begin{equation}  \int_{\partial P} u \ d\sigma \leq C \ L_{A, \sigma}
    u = 2 C  \Area(P) . \end{equation}
    This immediately gives {\it a priori} $C^{1}$ bounds on the restriction of $u$ to compact subsets
    of $P$, which can be applied to bound all derivatives in the interior,
    as in \cite{kn:D2}. In \cite{kn:D3} we showed that the estimates could
    be extended up to the boundary provided that $u$ satisfies a condition,
    called there an \lq\lq
    M-condition''. We recall the definition. Let $p,q$ be points in $P$ and
    write
    $q= p+ d\ \nu$ where $\nu$ is a unit vector. Write $V(p,q)$ for the difference
    in the derivative of $u$ in the $\nu$ direction evaluated at $q$ and
    $p$. Then $u$ satisfies an $M$ condition if for any pair $p,q$ in $P$
    such that the points $p-d\ \nu$ and $q+d\ \nu$ also lie in $P$ we have
    $ V(p,q) \leq M $. 
    
    Roughly speaking, the point of the present paper is to show that
    the solutions of our equation satisfy an {\it a priori} $M$-condition.
    But the  detailed strategy is considerably more complicated and we outline
    it now.  
    \begin{enumerate}
    \item In Section 2, we obtain an {\it a priori} $L^{\infty}$ bound on
    the solutions. (Notice that an $M$-condition implies such a bound, by
    an easy argument.)
    \item In Section 3, we work with points $p$ in  the neighborhood of an  edge of the polygon
    but away from the vertices. We obtain an {\it a priori} bound on a quantity
    $D(p)$, which is essentially equivalent to the $M$-condition for pairs
    of points $p,q$ not close to the
    the vertices. This is the core of  the paper: we use a \lq\lq blow-up''
    argument hinging on the compactness properties of  sets of bounded convex functions. 
    \item In Section 4 we use the arguments of \cite{kn:D3} to obtain bounds
    on the curvature tensor of our solution, away from the vertices. This
    requires some subsidiary arguments to control the Riemannian distance
    function, also using the bound on the quantity $D(p)$ from Section 3.
    \item In Section 5 we study the solutions in neighbourhoods of the vertices.
    First, using the  maximum principle and our $L^{\infty}$ bound from Section
    2, we obtain a two-sided bound on the volume form $\det u_{ij}$. We then
    use arguments similar to those in Section 3 to obtain an {\it a priori}
    bound on a quantity $E(t)$ which is related to the $M$-condition
    near the vertices.

    \end{enumerate}
    
    Many steps in this programme apply equally well to the general problem,
    with other functions $A$. The main obstacle in making this extension
    comes in the arguments of Section 3. Here we use a special estimate,
    for the case of constant $A$, proved in \cite{kn:D2}. If $u$ is any smooth
    convex function we can define a vector field $V$ by
    \begin{equation}     V^{i} = -\sum_{j} \frac{\partial u^{ij}}{\partial x_{j}}.\end{equation}
In \cite{kn:D3}, Theorem 2 we showed that, when the dimension is $2$ and the
function $A$ is constant, there is an {\it a priori} bound on the Euclidean
norm of $V$,
\begin{equation}   \vert V \vert_{{\rm Euc}} \leq C. \end{equation}
  Apart from the use of this estimate in Section 3, the only other restriction
  on $A$ comes when, in two places, we use the fact that $A$ is positive
  (in arguments involving the maximum principle) Thus the main result in this paper
extends to positive functions $A$ if one assumes an {\it a priori} estimate
(6). However the uses we make of (6) are to overcome some rather technical
difficulties which do not appear to be central to the problem, so there are
good grounds for hoping that the proofs may be extended to other functions
$A$ in time.

  Let us say a little more here about the central argument in Section 3 of
  the paper. The $M$-condition is essentially a local $C^{0}$ bound on a
  convex function. Suppose $u$ is a smooth convex function normalised at the origin. Then we can obviously choose a small scalar $\epsilon$ such that $\ufl=c u$ satisfies a fixed $C^{0}$ bound on a fixed ball about the origin
 and we have $\left(\ufl\right)^{ij}_{ij}= c^{-1} u_{ij}^{ij}$.
Then if we have a sequence $u^{(\alpha)}$ such that $\left(u^{(\alpha)}\right)^{ij}_{ij}$
is suitably small we can make a sequence of such re-scalings to get a new
sequence $\ufl^{(\alpha)}$, bounded in $C^{0}$ and with $\left(\ufl^{(\alpha)}\right)^{ij}_{ij}$
tending to zero. We can suppose the $\ufl^{(\alpha)}$ converge in $C^{0}$
over the interior of the ball and the question is: what can we say about
the limit? This is the basic idea used in Section 3, and also (more implicitly)
in Section 5. The complexities of the arguments, and the restrictions on
the results, are  related to our lack of understanding of this basic local
question---what
are the possible $C^{0}$ limits of solutions to the equation (1)?

    \

The main work of this paper finishes in Section 5. However, it is interesting
to complement the existence proofs with an understanding of what goes wrong
when the positivity hypothesis  is violated.  Section 6 is  a supplement
in which we discuss an explicit family of complete zero scalar curvature metrics, generalising
the Taub-NUT metric, and explain that these can be expected to arise as blow-up
 limits of solutions. We also discuss briefly the connection with \lq\lq
 collapsing'' phenomena in Riemannian geometry.

\section{The $L^{\infty}$ estimate}
   \subsection{An integral inequality}
   
   In this subsection we derive a geometric inequality for solutions of the
   equation (1). Consider the general case of a polytope $P\subset \bR^{n}$ with boundary measure $\sigma$ and a function $A$. We suppose that
$u$ is a solution to (1) satisfying Guillemin boundary conditions. Recall that
this implies that for any smooth test function $f$ on $\oP$ we have
$$  \int_{\partial P} f d\sigma = \int_{P} u^{ij} f_{ij} + A f d\mu . $$
Suppose that $\uu$ is a weakly-convex smooth function on $\oP$ and that
$\uu=u$ on an open subset $X\subset \oP$. Take $f= u-\uu$, so $f_{ij}$ vanishes
on $X$. Outside $X$ we have
$$   u^{ij} f_{ij} = u^{ij}(u_{ij}- \uu_{ij}) \leq u^{ij}u_{ij} = n, $$
since $\uu$ is convex. Thus
\begin{equation}  \int_{\partial P} u-\uu \ d\sigma \leq n \Vol(P\setminus X)+ \int_{P} A (u-\uu) \  d\mu .\end{equation}
Now start with an open set $X\subset\oP$, with piecewise-smooth boundary,
say. We define a convex function $\uu_{X}$ by
$$  \uu_{X} = \max_{p\in X\cap P} \lambda_{p} , $$ where for each $p\in P$
we write $\lambda_{p}$ for the affine-linear function defining the supporting
hyperplane of the graph of $u$. Thus $\uu_{X}=u$ on $X$, by convexity, and
an alternative definition is that $\uu_{X}$ is the least convex function
which restricts to $u$ on $X$.
We claim that the inequality (6) holds, with $\uu=\uu_{X}$. This is immediate if $\uu_{X}$ is smooth.
In the general case, we introduce a small parameter $\epsilon$ and let
$\uu_{X, \epsilon}$ be the standard mollification of $\uu_{X}$, using convolution
with a bump-function supported in the $\epsilon$- ball. (Notice that $\uu_{X}$ is defined
on all of $\bR^{n}$.) Then $\uu_{X,\epsilon}$ is again convex and is equal
to $u$ on the set $X_{\epsilon}\subset X$, defined by removing the $\epsilon$-neighbourhood
of the boundary of $X$. Then we can apply (6) to $\uu_{X,\epsilon}$, so 
    
    $$   \int_{\partial P} u-\uu_{X,\epsilon} d\sigma \leq n \Vol (P\setminus
    X_{\epsilon}) + \int_{P} A(u-\uu_{X,\epsilon}) \  d\mu. $$
    We take the limit as $\epsilon $ tends to $0$. The functions $\uu_{X,\epsilon}$
    converge uniformly to $\uu_{X}$ and the volume of $X_{\epsilon}$ tends
    to the volume of $X$ by our regularity assumption on the boundary of
    $X$. So, in sum, we have derived the geometric inequality
   \begin{equation} \int_{P} u - \uu_{X} \ d\mu\leq n \Vol(P\setminus X) + \int_{P} A
    (u-\uu_{X})\ d\mu. \end{equation}
    
    \

    (Notice that this inequality makes sense for arbitrary convex functions
    $u$. We can think of it as a partial \lq\lq weak form'' of the equation
    (1)
    , and the boundary conditions.)

\subsection{$L^{\infty}$ estimate: the main idea}

In this subsection and the next we will apply the ideas above to derive an {\it a priori}
bound for ${\rm max}\ \tu$, when $\tu$ is a normalised solution of (1), and
the boundary conditions, assuming
an $L^{1}$ bound on the restriction to the boundary.  Of course this maximum is attained at one of the vertices, so we fix
a vertex $q$ and seek to bound $\tu(q)$.  We can choose coordinates so that
$q$ is the origin and, near to  $q$, the polygon $P$ is the first quadrant
$\{x_{1}>0,x_{2}>0\}$. We write $(2 l_{1},0)$ and $(0,2 l_{2})$ for the coordinates
of the two vertices adjacent to $q$. We can arrange that the boundary measures
on these two edges are standard. Recall that
the $L^{1}$ bound on the boundary values of $\tu$ gives, by very elementary
arguments, bounds on $\tu$ and its first derivative in the interior of each
edge of the boundary. Let $u$ be the unique function obtained by adding an affine-linear
function to $\tu$ such that
\begin{itemize}
\item $\frac{\partial u}{\partial x_{1}}=-1$ at the midpoint $(l_{1},0)$
and $\frac{\partial u}{\partial x_{2}}=-1$ at $(0,l_{2})$.
\item The minimum value of $u$ on $P$ is $0$.
\end{itemize}
All of these preliminaries are   just to provide a convenient setting for
the main arguments.
Clearly, the bounds on the derivatives of $u$ at the midpoints mean that
it suffices to obtain an {\it a priori} bound on $u(0,0)$. 

The goal of this section is to prove
\begin{thm}
There is an $H$ depending only on $l_{1}, l_{2}, \Vert A \Vert_{L^{\infty}}$
and the integral of $u$ over the boundary of $P$ such that $u(0,0)\leq H$.
\end{thm}

\

To prove the Theorem, we will apply our inequality (8) to a $1$-parameter family of domains $X(h)$
in $P$. Define a function $\phi$ on $P$ by
$$  \phi = u-x_{1} \frac{\partial u}{\partial x_{1}}-x_{2}\frac{\partial
u}{\partial x_{2}}, $$
and set
$$ X(h)= \{ \ux: \phi(\ux)< h\}. $$
Thus $X(h)$ is the largest subset with the property that $\uu_{X(h)}(0,0)\leq
h$. As the parameter $h$ increases the domain $X(h)$ grows and once $h\geq
u(0,0)$ we have $X(h)=P$.  We write $\Omega(h)$ for the complement $P\setminus
X(h)$. We also write $\oh=u(0,0)$  and work with values $h<\oh$. Then the closure of $\Omega(h)$ meets the axes in a pair of line segments, from
the origin to $(\xi_{1}(h), 0), (0,\xi_{2}(h))$ respectively, say. 
We also write  $\uh_{1}= u(l_{1},0)+ l_{1}$ and $\uh_{2}=u(0,l_{2})+l_{2}$
and set
$\uh={\rm max} (\uh_{1},\uh_{2})$. Then if $h>\uh$ we have $\xi_{i}(h)\leq l_{i}$.
 These definitions are illustrated in the diagram.

 
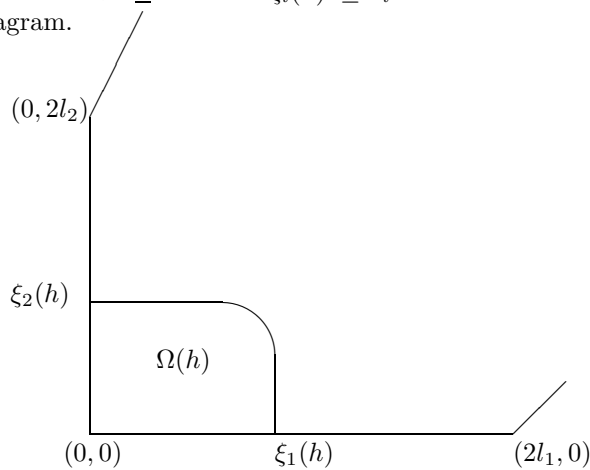
\begin{picture}(180,180)(-30,-30)
\put(0,0){\line(1,0){160}}
\put(0,0){\line(0,1){120}}
\put(160,0){\line(1,1){20}}
\put(0,120){\line(1,2){20}}
\put(-10,-10){$(0,0)$}
\put(-30,120){$(0,2l_{2})$}
\put(160,-10){$(2l_{1},0)$}
\put(70,-10){$\xi_{1}(h)$}
\put(-30,50){$\xi_{2}(h)$}
\put(25,25){$\Omega(h)$}
\put(0,0){\oval(140,100)[tr]}
\end{picture}

Let $\tau_{1,h}(t)$ be the affine-linear function
of one variable whose graph is the supporting hyperplane of the restriction
of $u$ to the $x_{1}$-edge at the point $\xi_{1}(h)$. Thus, by definition,
$  \tau_{1,h}(0)= h$. It follows from the definition that the restriction
of the function $\uu_{X(h)}$ to the axis is supported on the interval $[0,\xi_{1}(h))$
 on which it is equal to the affine-linear function
$\tau_{1,h}$.  Let
$$   G_{1}(h)= \int_{0}^{\xi_{1}(h)} u(t,0)-\tau_{1,h}(t) \ dt, $$
  and define $\xi_{2}(h)$ and $G_{2}(h)$ similarly.
  Then 
  $$  \int_{\partial P}u- \uu_{X(h)} d\sigma = G_{1}(h)+G_{2}(h). $$

We can now explain the main idea of our proof.  To begin with let us suppose that for $h\geq \uh$ the support of the function
$A$
does not meet $\Omega_{h}$. Our basic inequality (8) becomes
$$   G_{1}(h)+G_{2}(h) \leq {\rm Area} (\Omega_{h}). $$  
Elementary calculus gives the identities
\begin{equation} \frac{d G_{i}}{d h}= -\frac{1}{2} \xi_{i}(h). \end{equation}
In the standard model, where $u=x_{1} \log x_{1} + x_{2} \log x_{2} + {\rm
constant}$ say, it is easy to check that $\Omega_{h}$ is exactly the triangle
with vertices $(0,0), (\xi_{1}(h), 0), (0, \xi_{2}(h))$, so in this case
the area of $\Omega_{h}$ is $\xi_{1}\xi_{2}/2$. Suppose that, in our general
situation, we were able to show that $\Omega_{h}$ is not too different from
this triangle, in  that we have an inequality
\begin{equation}  {\rm Area}(\Omega_{h}) \leq \kappa \ \xi_{h} \xi_{2}(h), \end{equation}
for some fixed $\kappa$. For example, if we knew that $\Omega_{h}$ is contained
in the rectangle with vertices $(0,0), (\xi_{1}(h), 0), (0, \xi_{2}(h)),
(\xi_{1}(h), \xi_{2}(h))$ we could take $\kappa=1$. Under this supposition
we have
$$   G_{1}+G_{2} \leq \kappa \ \xi_{1} \xi_{2} \leq \kappa \left(
\frac{d G_{1}}{d h} + \frac{d G_{2}}{d h}\right)^{2},
$$
where we have used (9).
So the positive, decreasing, function $\Gamma=G_{1}+G_{2}$ satisfies the differential inequality

$$   \frac{d\Gamma}{dh} \leq - \sqrt{ \frac{ \Gamma}{\kappa}} $$
in the interval $\uh< h< \oh$. This gives
$$  \sqrt{\Gamma}(\uh)- \sqrt{\Gamma}(h) \geq \frac{1}{\sqrt{\kappa}} (h-\uh), $$
and thus, since $\Gamma$ tends to $0$ as $h$ tends to $\oh$,
$$\oh-\uh\leq \sqrt{ \kappa \Gamma(\uh)}. $$
Since $\Gamma (\uh)$ is dominated by the integral of $u$ over the boundary of
$P$ this gives the desired bound on $\uh=u(0,0)$.

\

To turn this idea into a complete proof we need to overcome two difficulties.
The first is to incorporate the term involving the function $A$ in our basic
inequality. This is relatively easy. The second, more fundamental, difficulty
is that the author does not know how to obtain a universal inequality of
the form (10) that we used above, although it seems very reasonable to expect
this to be true. Thus the actual proof, which we give in the
next section, is more complicated since it is based on  a weaker assertion
than (10) (Lemma 2 below).

\subsection{The detailed proof}
We begin with  some elementary calculus associated
to a convex function of one variable. This will be applied to
the boundary values of our function $u$, but to simplify notation consider
first a strictly convex, smooth,  function $U(t)$ on an interval $(0,2l)$
with $U'(l)=-1$. For $h\geq U(l)+ l$ we define $\xi(h)\in (0,l)$ as above,
i.e. so that the affine-linear function $\tau_{h}$ whose graph is tangent
to the graph of $U$ at $\xi$ has $\tau_{h}(0)=h$.
Let $D(h)$ be the point where the affine-linear function $\tau_{h}$
vanishes. In other words, the line joining the two points $(h,0)$ and
$(0,D(h))$ is tangent to the graph of $U$ at the point $( \xi(h), U(\xi(h))\
)$. Let $z(h)= D(h)/h$, so $z(h)^{-1}=- u'(\xi(h))$.
\begin{lem}
In this situation
$$   \xi(h)= \frac{D^{2}}{D-hD'}. $$ 
\end{lem}
This is a calculus exercise for the reader.

 \

  Now return to our function $u$ of two variables. We extend the notation
  above in the obvious way, so we have functions $D_{i}(h), z_{i}(h)$ defined
  for $\uh \leq h \leq \oh$. For these values of $h$, we let $\Delta_{h}$
  be the triangle in the $(x_{1}, x_{2})$ plane with vertices $(0,0), (D_{1}(h),0), (0, D_{2}(h))$. 
\begin{lem}
For any $h\in (\uh,\oh)$ we have $\Omega(h)\subset \Delta_{h}$.
\end{lem}
To see this, consider a point $p$ in $\Omega(h)$. Let $\pi$ be the affine-linear
function defining the supporting
hyperplane of  $u$ at the point $p$ and write $h^{*}=\pi(0,0)$. The condition that
$p$ lies in $\Omega_{h}$ is the same as saying that $h^{*}\geq  h$. Consider
the restriction of $\pi$ to the $x_{1}$-axis. By convexity we have $\pi(t,0)\geq u(t,0)$ for all $t$, in particular $h^{*}\leq u(0,0)$ and so $\xi(h^{*})$
is defined. Then $\tau_{1,h^{*}}$ and the restriction of $\pi$
are two affine-linear functions of one variable, equal at the origin.
We must have $\pi(t,0)\geq \tau_{1,h^{*}}(t) $ for all $t>0$, for otherwise
$\pi(t,0)\leq \tau_{h^{*}}(t)$ for all $t>0$, which is a contradiction when $t=\xi_{1}(h^{*})$. Thus $\pi(D^{*}_{1},0)=0$ for some $D^{*}_{1}<D_{1}(h)$.
Similarly $\pi(0,D^{*}_{2})=0$ for some $D^{*}_{2}<D_{2}(h)$. Now $\pi(0,0)>h>0$ so
the region in $P$ where $\pi>0$ is the triangle $\Delta^{*}$ with vertices
$(0,0), (D^{*}_{1},0), (0, D^{*}_{2})$. At the original point $p$ we have $u(p)=\pi(p)$. Since
$u$ was normalised so that $u\geq 0$ we have $\pi(p)\geq 0$ and so $p$
lies in $\Delta^{*}$. But $\Delta^{*}$ is contained in $\Delta_{h}$, since
$D^{*}_{i}< D_{i}$, so $p$ lies in $\Delta_{h}$, as required.

\

Next we look at  the term involving the function $A$. For
$h\geq \uh$ we write $f_{h}=u-\uu_{X(h)}$. So $f_{h}$ is a positive
function, supported in
the set $\Omega_{h}$. We have to estimate the term
$$  \int_{\Omega_{h}} A f_{h} , $$
appearing in the inequality (7). Set $\alpha= \Vert A \Vert_{L^{\infty}}$ so
\begin{equation}  \int_{\Omega_{h}} A f_{h} \leq \alpha J(h), \end{equation}
where 
\begin{equation}  J(h)= \int_{\Omega_{h}} f_{h}. \end{equation}

\begin{lem}
With notation as above,
$$  \frac{d J(h)}{dh}= - \frac{1}{3} {\rm Area}\ ( \Omega_{h} ). $$
\end{lem}

This is the two-dimensional analogue of the elementary identity (9).  
To prove it we work in polar coordinates, writing $u(r,\theta)$. Consider a ray through the origin, with fixed $\theta$. The restriction of $u$ to
this ray is a convex function of $r$ and there is a unique point $r=R(\theta)$
where $u- r\frac{\partial u}{\partial r}= h$. From the definitions, this
is a point on the boundary  $\partial \Omega_{h}$ and the restriction of
the function  $\uu_{X(h)}$ to the intersection of $\Omega_{h}$ and this ray is
$$\uu_{X(h)}(r,\theta) = h+\frac{r}{R} (u(R,\theta)-h). $$
It follows that
$$  \frac{\partial}{\partial h} \uu_{X(h)}(r,\theta)= 1- \frac{r}{R}. $$
Thus
$$  \frac{d}{dh} J(h)= \int_{\theta=0}^{\theta=\pi/4} \int_{r=0}^{r=R(\theta)} (1-\frac{r}{R})\  r dr d\theta. $$
Performing the $r$ integral this is
$$ \frac{d}{dh} J = \frac{1}{6} \int_{\theta=0}^{\theta=\pi/4} R(\theta)^{2} d\theta, $$
while the area of $\Omega_{h}$ is given by the usual formula
$$  {\rm Area}\ (\Omega_{h}) = \frac{1}{2}\int R(\theta)^{2} d\theta. $$

\

Combining the two lemmas above, we get
\begin{cor}
$$J(h) \leq \frac{1}{6} \int_{h}^{\oh} D_{1}(h) D_{2}(h) \ dh. $$
\end{cor}
This follows immediately from the co-area formula and the  facts that the area of $\Delta_{h}$ 
is $D_{1}D_{2}/2$ and that $J(h)\rightarrow 0$ as $h\rightarrow \uh$.

\

Now define $I_{1}(h)$, for $\uh\leq h < \oh$ by
\begin{equation}  I_{1}(h) = \int_{h}^{\oh} D_{1}(h)^{2}\  dh,  \end{equation}
 and define $\lambda_{1}(h)$ by the equation
\begin{equation}   G_{1}(h)= \frac{\lambda_{1}(h)}{2} D_{1}(h)^{2} + \frac{\alpha}{12}
I(h). \end{equation}
(The reason for the choice of factor $\alpha/12$ will appear shortly.)

Define $I_{2}$ and $\lambda_{2}$ similarly.
\begin{lem}
$\lambda_{i}(h) \rightarrow 0$ as $h\rightarrow \oh$
\end{lem}
This is straightforward to check, using the known behaviour of $u$ at the
origin. We omit the details.

\

Now we can proceed to the core of the  proof, which has two parts. The
first is stated in 
\begin{prop}
Suppose $\uh\leq h \leq\oh$ and $\lambda_{1}(h), \lambda_{2} (h)$ are both
positive. Then
$$   \lambda_{1}(h) \lambda_{2}(h) \leq 1. $$
\end{prop}

 Using Lemma 2, (11) and Corollary 2, our basic inequality (8)
gives, for any $h\in (\uh,\oh)$,
$$    G_{1}(h)+ G_{2}(h) \leq \frac{1}{2} D_{1}(h) D_{2}(h) + \frac{\alpha}{6} \int_{h}^{\oh}
D_{1}(s) D_{2}(s) \ ds. $$
Using the inequality $D_{1} D_{2} \leq \frac{D_{1}^{2}+D_{2}^{2}}{2}$ and
the definition of $I_{i}(h)$ we
get
\begin{equation}  G_{1}(h)+G_{2}(h) \leq \frac{1}{2} D_{1}(h) D_{2}(h)+ \frac{\alpha}{12}
(I_{1}(h)+I_{2}(h)). \end{equation}
So we have, from the equations defining $\lambda_{i}$,
$$  \frac{\lambda_{1}}{2} D_{1}^{2} +\frac{ \lambda_{2}}{2}D_{2}^{2} \leq
D_{1} D_{2}, $$
for each $h\in (\uh,\oh)$. In other words
$$  \frac{1}{2} \left( \lambda_{1} \frac{D_{1}}{D_{2}} + \lambda_{2} \frac{D_{2}}{D_{1}}
\right) \leq 1. $$
From the arithmetic-geometric mean inequality we see that if $\lambda_{1}(h),
\lambda_{2}(h)$ are both positive then
$\lambda_{1}\lambda_{2}\leq 1$.  

\

For the second part, we derive differential equations involving  the functions
$\lambda_{i}(h), z_{i}(h)$. For clarity we suppress the suffix $i$ temporarily, and
denote derivatives with respect to $h$ by a prime symbol. Write $c=\alpha/12$ and recall that
$  G=\frac{\lambda}{2} D^{2} + c I $, where $\frac{dI}{dh}= -
D^{2}$. Thus, using (9), 
$$   -\frac{\xi}{2}= G'= \lambda D D'+ \frac{1}{2}\lambda' D^{2} - c D^{2}.
$$
By Lemma 1 this gives,
$$  \frac{ D^{2}}{2(h D'-D)}= \lambda D D' + \frac{1}{2} \lambda' D^{2} - cD^{2}.
$$
Since $D=z h$, we have 
$$   hD'-D=h^{2} z'<0,$$
and our equation becomes
$$   \frac{z^{2}}{2z'}= \lambda z h ( h z'+z) + \frac{\lambda'}{2} z^{2} h^{2}-
cz^{2} h^{2}. $$
This leads to
$$   \frac{z}{2h^{2}} = z'\left(\lambda z'+ z(\frac{\lambda}{h}- c+\frac{\lambda'}{2})
\right). $$
Recall that $z>0$ and $z'<0$. Suppose $K$ is any fixed positive number. For
any $A,B>0$ 
we have $K\sqrt{AB}\leq \frac{1}{2}( K^{2} A +  B)$.
We apply this to the right hand side of (15), with $A=-z'$ and $B=-(\lambda
z' + z(\lambda/h - c+\lambda'/2)$. We deduce that
\begin{equation}  z'(K^{2}+\lambda)+ z( \frac{\lambda}{h}-c +\frac{\lambda'}{2}) \leq - \sqrt{2} K \frac{\sqrt{z}}{h}.
\end{equation}

\

\begin{prop}
Suppose $z(h), \lambda(h)$ are functions defined on an interval $(h_{0},
\oh)$ with the following properties
\begin{enumerate}
\item $z(h)>0$ and $z'(h)<0$ for all $h$.
\item $z,\lambda$ satisfy the differential inequality (16) above, for some
$c>0$ and all $ K>0$.
\item For some $C>0$, and all $h$, we have $$z\leq C/h^{2}.$$ 

\item $z(h)$ and $\lambda(h)$ tend to $ 0$ as $h\rightarrow \oh$.
\end{enumerate}
Write $b=2\sqrt{2}-1$. If we fix any $K>2c\sqrt{C} $ and if we set $h_{1}={\rm max}(h_{0},
\frac{3 K\sqrt{C}}{b})$, then we have $K^{2}+\lambda(h)>0$ for all $h\geq h_{1}$ and
$$  \int_{h_{1}}^{\uh} \frac{1}{(K^{2}+\lambda)^{3/4}} \ \frac{dh}{h} \leq \frac{12}{b K}
\ z(h_{1})^{1/2}( K^{2} + \lambda(h_{1}))^{1/4} . $$ 
\end{prop}

We fix $K$ as stated. Multiplying the inequality (16) by $2z$, we have
$$  2 z z'(K^{2}+\lambda)+ \lambda'z^{2} + 2 z^{2}\left( \frac{\lambda}{h}-c\right)
\leq- \frac{2 \sqrt{2} K}{h} z^{3/2}. $$
This is
$$  \frac{d}{dh} \left( z^{2} (K^{2}+\lambda)\right)\leq - \frac{2 \sqrt{2}K}{h} z^{3/2}+
2z^{2} \left( c-\frac{\lambda}{h}\right). $$
Now the inequalities $z\leq C h^{-2}$ and $K>2c\sqrt{C}$ imply that
$2z^{2}c < \frac{K}{h} z^{3/2}$ so we have
\begin{equation} \frac{d}{dh}\left( z^{2}(K^{2}+\lambda)\right) \leq -\frac{b K}{h} z^{3/2}
- \frac{2\lambda}{h} z^{2}. \end{equation}

Write $F=K^{2} +\lambda$. We want to show that $F(h)$ is positive for $h\geq
h_{1}$, so  we suppose that $F(h_{2})< 0$ for some $h_{2}>h_{1}$ and seek
a contradiction.
Since $F\rightarrow K^{2}>0$ as $h\rightarrow \oh$ there is an $h_{3}\geq
h_{2}$
with $F(h_{3})=0$ and $F'(h_{3})\geq 0$. This means that
$$  \frac{d}{dh}(z^{2} F) = z^{2}F'+ (z^{2})'F$$ is positive when $h= h_{3}$.
 So, with $\lambda=\lambda(h_{3}), z=z(h_{3})$ we have
$$  2\frac{\lambda}{h}  z^{2}\leq  -\frac{bK}{h} z^{3/2} .$$
But $\lambda=-K^{2}$ (since $F=0$) and we get
$$   K z^{1/2}\geq \frac{b}{2} . $$
Since $z\leq C/h_{3}^{2}$, we have
$$  \frac{\sqrt{C}}{h_{3}} \geq \frac{b}{2K} . $$
But this contradicts the assumption that $h_{1}\geq 3K\sqrt{C}/ b$, since
$h_{3}\geq h_{1}$. So we have established that
$F(h)>0$ for $h>h_{1}$.

\

Now if $h>h_{1}$ we  have
$$  -\frac{\lambda}{h}\leq \frac{K^{2}}{h}\leq \frac{K^{2}}{h_{1}} $$
and so
$$   -2 z^{2} \frac{\lambda}{h} \leq \frac{2\sqrt{C}K^{2}}{h_{1}}\ \frac{ z^{3/2}}{h}. $$
So we obtain from (17) that
$$  \frac{d}{dh}( z^{2} F) \leq (\frac{2K^{2}\sqrt{C}}{h_{1}}- bK) \frac{z^{3/2}}{h}
. $$
By the choice of $h_{1}$, this gives,
$$  \frac{d}{dh}(z^{2} F) \leq - \frac{bK}{3} \frac{z^{3/2}}{h}. $$
 Write $w=z^{2}F$ so
the above inequality is
$$  \frac{dw}{dh}\leq - \frac{bK}{3 h} \frac{w^{3/4}}{F^{3/4}}. $$
That is
   $$\frac {d}{dh} w^{1/4}\leq -\frac{bK}{12} F^{-3/4} \frac{1}{h}. $$
   
 We know that $w(h)$ tends to zero
   at $\uh$ so we can integrate this with respect to $h$ to obtain the inequality
   stated in Proposition 2.
   
   \
   
   Propositions 1 and 2 are the essential parts of the proof of Theorem 2,
    and it remains
   now to put together the various components. Of course we want to apply Proposition 2 to
   the functions $\lambda_{i}, z_{i}$ associated to the two edges, with $h_{0}\geq
   \uh$. We have
   to show that there is a bound $z_{i}\leq Ch^{-2}$. Consider
  the quantity $z_{1}h^{2}/2=h D_{1}(h)/2$. This is the  integral of the affine-linear
  function
  $\tau_{1,h}$ from $0$ to $D_{1}(h)$. The value $u(l_{1}, 0)$ is controlled
  by the integral of $u$ over the boundary so there is no loss in supposing
  that $h_{0}\geq 2 u(l_{1}, 0)$. This implies that, for $h\geq h_{1}$ we
  have
  $D_{1}(h)\leq 2l_{1}$ and, since $\tau_{1,h}(t)\leq u(t,0)$ for $t$ in the interval $[0,
  D_{1}(h)]$, we   get
  $$ \frac{1}{2}  h D_{1}(h) \leq \int_{\partial P} u. $$
  This gives the desired bound $z_{1} \leq C h^{-2}$, and similarly of course
  for $z_{2}$.
   Thus we may apply Proposition $2$, with
  a suitable fixed $K$ determined by $C$ and $c= \Vert A \Vert_{L^{\infty}}/12$. Now consider the functions
  $w_{i}(h)= z_{i}^{2}(K^{2}+\lambda_{i}(h))$. By construction $z_{i}(h) \leq
  1$ for $h\geq \uh$ so
  $w_{i}\leq K^{2} + z_{i}^{2} \lambda_{i}=K^{2}+ \lambda_{i} D_{i}^{2}/h^{2}$.
  By the definition of $\lambda_{i}$ we have
  $$  \lambda_{i} D_{i}^{2}= 2(G_{i}- cI_{i})\leq 2 G_{i}. $$
  Obviously the functions $G_{i}$ are bounded by the integral of $u$ over $\partial P$.
  So we obtain
  $$  w_{i} \leq K^{2} + \frac{2}{h^{2}}\int_{\partial P} u. $$
  This gives an upper bound on $w_{i}(h_{1})$, since $h_{1}$ is determined
  by $C$ and $c$.
  In sum, Proposition 2 tells us that there are $h_{1}, L, K$, all determined
  by the integral of $u$ over the boundary, such that $K^{2}+ \lambda_{i}(h)>0$ if $h>h_{1}$ and
  $$  \int_{h_{1}}^{\oh} \frac{1}{(K^{2} + \lambda_{i})^{3/4}} \frac{dh}{h}
  \leq L. $$
  Change variable by writing $h=e^{t}$ and let $\oh=e^{\ot}, h_{1}=e^{t_{1}}$.
  Then we have
  $$ \int_{t_{1}}^{\ot} \frac{1}{(K^{2}+\lambda_{i})^{3/4}}\  dt \leq L. $$
  So the measure of the set in $[t_{1}, \ot]$ where $\lambda_{i}\leq 1$ is
  at most $L(K^{2}+1)^{3/4}$. Thus if $\ot$ were bigger than $t_{1}+2L(K^{2}+1)^{3/4}$
  there would have to be a point where $\lambda_{1}>1$ and $\lambda_{2}>1$.
  But this would contradict Proposition 1. So we conclude that
  $\ot\leq t_{1}+ 2L(K^{2}+1)^{3/4}$ or in other words
  $$   \oh \leq h_{1} \exp(2L(K^{2}+1)^{3/4}), $$
  and we have  proved Theorem 2.

\section{Edges}
\subsection{Preliminaries}
We now come to the central topic of this paper. Consider a symplectic potential
$u$ on a polygon $P$ and let $E$ be an edge of 
$P$. Choose an outward-pointing  vector $\nu$ transverse to $E$---say the Euclidean
normal. Suppose that $p$ is a point of $P$ such that the ray $\{ p + t\nu:
t>0\}$ meets the edge $E$ in a point $q=p+s \nu$, for $s=s(p)$. Let $\lambda_{p}$ be the
affine-linear function defining the supporting hyperplane to $u$ at $p$;
so the difference $u-\lambda_{p}$ vanishes to first order at $p$. Then we
define
\begin{equation}  D(p)= \frac{u(q)-\lambda_{p}(q)}{s(p)}. \end{equation}
The goal of this section is to obtain an {\it a priori} bound on $D(p)$,
under mild hypotheses. It is easy and elementary to go from this to an
\lq\lq M-condition'' formulation, as we will explain in Section 4. When we
want to indicate the dependence on the function $u$ we write $D(u;p)$.

We will want to have this {\it a priori} bound in the context of the continuity method of \cite{kn:D3};
when we have a sequence $(P^{(\alpha)}, A^{(\alpha)}, \sigma^{(\alpha)})$
of data sets and convex functions $u^{(\alpha)}$. In this context the data
sets will converge in the obvious sense as $\alpha\rightarrow \infty$. To
simplify notation we will often omit the index $\alpha$, and just write
$P,A,\sigma$, where it is clear that the quantities involved (for example
the diameter of $P^{(\alpha)}$) satisfy a uniform bound in the sequence.
We suppose that $u^{(\alpha)}$ is normalised (in the sense of Section
1) at the centre of mass of $P^{(\alpha)}$.
The main result we prove is
\begin{thm} Suppose that the data sets $(P^{(\alpha)}, A^{(\alpha)}, \sigma^{(\alpha)})$
converge as $\alpha\rightarrow \infty$ and that the sequence $u^{(\alpha)}$ satisfies uniform bounds
 $${\rm max}_{P} u^{(\alpha)} \leq C_{0}.$$
 $$  \vert V^{(\alpha)}\vert_{{\rm Euc}} \leq C_{1}. $$
Fix any $\delta>0$. There are $D_{\delta}, s_{\delta}$ with the following
property. If $p$ is
a point in $P^{(\alpha)}$ with $s(p)\leq s_{\delta}$ and the distance of
 $ p+s(p) \nu$ from the vertices of $P^{(\alpha)}$ is greater than $\delta$
 then $D(p) \leq D_{\delta}$.
 \end{thm}
 
 Here $V^{(\alpha)}$ is the vector field associated to $u^{(\alpha)}$ by the
 formula (5). Strictly speaking we should write $s_{\alpha}(p)$ etc., since
 these quantities depend on $P^{(\alpha)}$, but we hope that the meaning is
 clear. The arguments in this section do not depend strongly on the  $L^{\infty}$ bound on $u$, as opposed to a $L^{1}$ bound on the boundary
 value. The former is only used once, in the proof of Proposition 3, and could
 be avoided with a little extra work. Of course, we showed in the previous
 Section that the two conditions are in fact equivalent.

 We will now explain the main idea of the proof of Theorem 3. We suppose,
 on the contrary, that there is a sequence of points $p_{\alpha}$ and $D(p_{\alpha})=D_{\alpha}$ tends to infinity. Let us also suppose that $p_{\alpha}$ is the \lq\lq worst''
point, maximising the function $D$ for each fixed $\alpha$ and that the sequence $p_{\alpha}$ stays a definite distance from the vertices. Then by performing
a sequence of affine transformations, adding suitable affine linear functions
and multiplying by $D_{\alpha}^{-1}$ we can obtain a sequence of convex functions
$\ufl^{(\alpha)}$ defined on large convex subsets of a half-plane $\{ x_{1}\geq
0\}$ with $\ufl^{(\alpha)}(0,0)=1$ and $\ufl^{(\alpha)}$ attaining its minimum
$0$ at the point $(1,0)$. Our overall strategy is to obtain a contradiction by showing
that the $\ufl^{(\alpha)}$ have a $C^{0}$ limit and making various arguments
with this, using the fact that $p_{\alpha}$ is the \lq\lq worst'' point.
 Two of the issues we have to deal with are
\begin{itemize}
\item In reality we need to use a more complicated definition of \lq\lq worst''
point, because of the constraint involving $\delta$.
\item We have to contend with the affine invariance of the problem, in choosing
the affine transformations to define $\ufl^{(\alpha)}$ appropriately. This
is the choice of the parameter $\lambda$ below.

\end{itemize}

The vector field $V$ associated to the convex function $u$  comes in to our arguments at a number of places and
we will now recall two relevant points of theory.   The first is that the vector
field encodes the boundary conditions. Expressed in terms of
coordinates, this says that the normal component of $V$ at a point of an edge
is fixed by the given measure $\sigma$. The second is that if we write
$L=\log \det u_{ij}$ and introduce Legendre transform coordinates
$\xi_{i}= \frac{\partial u}{\partial x_{i}}$ then
\begin{equation}  V^{i}= \frac{\partial L}{\partial \xi_{i}} . \end{equation}
This leads to a basic principle which will be important in our arguments.
Suppose we have a bound on one component of the vector field: $\vert V^{2}\vert \leq C$ say. Then if, over a portion $\Gamma$ of  a contour $\{ \xi_{1}= {\rm constant}\}$
the partial derivative $\frac{\partial u}{\partial x_{2}}$ varies by a bounded
amount $b$ say, then the ratio 
$$  \frac{\det u_{ij}(\gamma)}{\det u_{ij}(\gamma')}$$
is bounded by $ e^{C b}$ for any $\gamma, \gamma'\in \Gamma$.

\

It is useful to have in mind a standard model for the boundary behaviour
given by the function
$$  u_{0}(x_{1}, x_{2})= x_{1} \log x_{1} + x_{2}^{2}, $$
on the half-plane $\{ x_{1}>0\}$. The associated vector field $V$ has components
$V^{1}=
1, V^{2}=0$, and the function satisfies (1) with $A=0$.  Then in this case
$D(p)$
is equal to $1$ for all $p$. Now apply an affine transformation and, for $a\in \bR$, set
$$  u_{a}(x_{1}, x_{2})= x_{1} \log x_{1} + (x_{2}- a x_{1})^{2}. $$
This also satisfies (1) with $A=0$, and the same boundary conditions.
At the point $p=(1,0)$ we have $D(p)= a^{2}+1$. Since $a$ can be made arbitrarily
large, this shows that there is no way to derive an {\it a priori} estimate
for $D(p)$ using only \lq\lq local'' information. However in this example
we have $V^{1}=1, V^{2}= a$, so the parameter $a$ is detected by the tangential
component of the vector field. This may help in  understanding our main proof
below, which uses the bound on $V$---obtained in our application by global arguments--to help
control the quantity $D(p)$.

\

With this outline of the strategy in place, we now proceed in more detail.
We suppose the polygon $P$ has an edge given by 
$$ E= \{(x_{1}, x_{2}): x_{1}=0, a\leq x_{2} \leq b\}, $$
and that $P$ lies in the half-plane $\{ x_{1} >0\}$. We can suppose that
the measure $d\sigma$ on this edge is the standard Lebesgue measure $dx_{2}$. We consider a solution
$u$ of our equation (1) , where $\vert A \vert\leq C_{2}$. We fix a $\delta>0$
and consider points
$   p= (s, t) $ with $s\leq s_{\delta}$ and $a+\delta \leq t \leq
b-\delta$.  We assume a bound
on the boundary integral, as in the statement of Theorem 3. It is then easy
to see that we can choose $s_{\delta}$ such that all such points $p$ lie in
$P$ and that  $ s D(p)$ satisfies a fixed bound
\begin{equation}  s D(p)\leq C_{3}. \end{equation}
Similarly we have  an elementary bound on the derivatives in the \lq\lq tangential
direction'': for any two points $p,p'$ satisfying the conditions above

    \begin{equation}  \vert \frac{\partial u}{\partial x_{2}}(p)-\frac{\partial
    u}{\partial x_{2}}(p')\leq C_{4} . \end{equation}

  Given $s,t$ as above, we define the convex function $u^{*}$ by normalising $u$
at $p=(s,t)$. That is, $u^{*}$ is given by adding an affine-linear function to
$u$ and $u^{*}$ attains its minimum value $0$ at $p$. Then, by definition, 
$$  D(p)=s^{-1} u^{*}(0,t). $$
 Recall that the Guillemin
boundary conditions around a vertex imply that the integral of the second
derivative $u_{22}=\frac{\partial^{2} u}{\partial_{x_{2}}^{2}}$,
evaluated along the edge, {\it diverges} at
each end point. Let $\lambda_{0}=\min( t-a, b-t)$ and for $\lambda<\lambda_{0}$
define
$$    I(\lambda)=\int_{t-\lambda}^{t+\lambda} u_{22}(\tau,0) d\tau. $$
Then $I(\lambda)$ is an increasing function, equal to $0$ when $\lambda=0$
and tending to infinity as $\lambda\rightarrow \lambda_{0}$. The same is
true of $\lambda I(\lambda)$ so there is a unique $\lambda$ such that
$$  \lambda I(\lambda) = s D(p)/2. $$
The reason for this choice will appear presently. 

\

The next proposition will be used to handle the potential \lq\lq end-point''
difficulties alluded to above.
\begin{prop}
There is a constant $c$, depending only on  $\delta, C_{0}, C_{2}, C_{3}, C_{4}, \Vert A \Vert_{L^{\infty}}$
and the geometry of the polygon $P$, such that $s(p) \leq c \lambda^{-2}$.
\end{prop} 
We need an elementary lemma.
\begin{lem}
There exists $\kappa>0$ with the following property. Let $f$ be any positive convex
function on $[-1,1]$ with $f(0)=1$ and with $\int_{-1}^{1} f''(t) dt =1$.
Let ${\cal C}_{f}$ be the set of affine-linear functions $\sigma$ such that
\begin{itemize}
\item
$\sigma(t)\leq f(t)$ for all $t\in [-1,1]$
\item {\it either} $\sigma\leq 0$ on $[-1/2, \infty)$ {\it or} $\sigma \leq 0$ on
$(-\infty, 1/2]$.
\end{itemize}
Define $g(t)= \sup_{\sigma \in {\cal C}_{f}} \lambda(t)$. Then
$$\int_{-1}^{1} f-g \ dt \geq \kappa. $$
\end{lem}

We leave the proof as an exercise.

\

To prove the Proposition we claim
first that if $D(p)$ is large then $s/\lambda$ is small. For if $\lambda\leq \delta/2$, then
$$I(\lambda)\leq \int_{t-\lambda}^{t+\lambda} u_{22}(0,\tau) d\tau\leq C_{4}.$$ So $s D(p)  = \lambda I(\lambda) \leq \lambda C_{4}$ and $s/\lambda \leq  C_{4}/S$. 

On the other hand if $\lambda>\delta/2$ then
$$   s/\lambda< 2s/\delta\leq 2C_{3}/D(p). $$

Now let $W$ be the wedge-shaped region
$$   W=\{ (x_{1}, x_{2}): \vert x_{2}-t \vert \leq \frac{\lambda}{2 s} (x_{1}-s),
$$
and let $X$ be the intersection of $P$ with $W$. It is clear that, when $s/\lambda$
is small the centre of mass of $P$ lies in $X$. 
 Since $X$ is convex, this implies that $\uu_{X}\geq 0$, where $\uu_{X}$
 is the function defined in Section 2.  So
$\vert u-u_{X}\vert \leq u$ and
\begin{equation} \int_{P\setminus X} (u-u_{X}) A  \leq \Vert A \Vert_{L^{\infty}} \max u\ \ 
\Area(P\setminus X). \end{equation}
Using our bound on $\max u$ we see that this integral is bounded by a multiple
of the area of $P\setminus X$. So (8) gives
$$ \int_{\partial P} u-u_{X} \leq c_{1} \Area(P\setminus X),$$
for some fixed $c_{1}$.
It is also clear that the area of $P\setminus X$ is bounded by a multiple
of $s/\lambda$, with the multiple depending only on the geometry of $P$. So we have
\begin{equation} \int_{\partial P} u-\uu_{X} \leq c_{2} \left( \frac{s}{\lambda}\right). \end{equation}
 
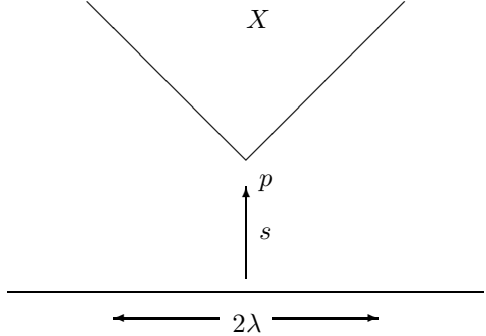
 \begin{picture}(180,180)(-30,-30)
\put(0,0){\line(1,0){180}}
\put(95,40){$p$}
\put(80,-10){\vector(-1,0){40}}
\put(100,-10){\vector(1,0){40}}
\put(85,-15){$2\lambda$}
\put(90,50){\line(1,1){60}}
\put(90,50){\line(-1,1){60}}
\put(90,100){$X$}
\put(90,5){\vector(0,1){35}}
\put(95,20){$s$}
\end{picture}

We now make a similar argument to that in Lemma 2. Let $q$ be a point of $X$, not equal to the centre of mass of $P$, and
let $\pi$ be the affine-linear function defining the supporting hyperplane
to $u$ at $q$.
Thus the zero set of $\pi$ is a line in the plane. Since $\pi(q)>0$
and $\pi(p)<0$ this line separates the points $q$ and $p$.
It follows that the restriction of $\pi$ to the $x_{2}$ axis is negative
on either the interval $(-\infty, -\lambda/2]$ or on the interval $[\lambda/2,
\infty)$. Lemma 5 (after suitable rescaling) implies that
$$  \int_{t-\lambda}^{t+\lambda} (u-\uu_{X})(0,\tau) d\tau \geq \kappa D(p) s \lambda. $$
 So from (23) we obtain
$$   \kappa D(p) \lambda \leq c_{2} \frac{s}{\lambda}, $$
which gives
$$     D(p) \leq \frac{c_{2}}{\kappa} \lambda^{-2}, $$
as required.

\

\

\

For $t$ in the interval $(a+\delta,b-\delta)$ define
$$  \Lambda(t)= \min (\vert (a+\delta)-t\vert, \vert (b-\delta)- t \vert) . $$
Set
$$  \mu= \max_{s,t} \Lambda(t)^{2} D(s,t), $$
where $t$ runs from $a+\delta$ to $b-\delta$ and $s$ runs from $0$ to $s_{\delta}$.
Choose a point $(s_{0}, t_{0})$ where the maximum is attained.
Of course, this is the concept of the \lq\lq worst point'' alluded to above.
 Write
$\lambda_{0}=\lambda(s_{0}, t_{0})$ and $D_{0}=D(s_{0}, t_{0})$. By the preceding Proposition we have
\begin{equation}    \lambda_{0} \leq c \sqrt{\mu} \Lambda(t_{0}). \end{equation}
In other words if, as we suppose, $\mu$ is large the \lq\lq scale'' $\lambda_{0}$
in the $x_{2}$ direction is small compared with the distance to the end points
$a+\delta, b-\delta$. Also for any fixed $R$ we have
\begin{equation} D(s',t')\leq D_{0} \frac{1}{1- R \sqrt{c}/{\sqrt \mu}} \end{equation}
for any $t'$ with $\vert t_{0}-t'\vert \leq R \lambda_{0}$.

\

\

We now define a new convex function $\ufl$ by rescaling. We set
\begin{equation}\ufl(x_{1}, x_{2}) = \frac{1}{D_{0} s_{0}} u^{*}( s_{0} x_{1}, t_{0}+\lambda_{0} x_{2}). \end{equation}
Thus $\ufl$ is defined on a polygon $\Pfl$ in the half-space $\{x_{1}>0\}$, which depends
on $P,s_{0}, t_{0},\lambda_{0}$. This polygon contains a rectangle $$Q=\{0<x_{1}< c s_{0}^{-1},
\vert x_{2}\vert < c \lambda_{0}^{-1}\}$$ for some fixed $c$ depending on
$P$ and $\delta$. Since we are supposing that $D_{0}$ is large, both $\lambda_{0}$
and $s_{0}$ are small (by (20) and Proposition 3). Thus the rectangle is
very  large.  By construction $\ufl$ attains its minimum value $0$ at
the point $(1,0)$ and $\ufl(0,0)=1$.  The choice of $\lambda$ transforms
into the condition that
\begin{equation}    \int_{-1}^{1}  (\ \ufl\ )_{22}(0,t) dt = 1/2. \end{equation}
The function $\ufl$ satisfies an equation
$$ \sum \frac{\partial^{2} \ufl^{ij}}{\partial x_{i}\partial x_{j}} = - \Afl, $$
in $\Pfl$,
where $\Afl(x_{1}, x_{2}) = D_{0} s_{0} A(s_{0} x_{1}, t_{0}+\lambda_{0} x_{2})$. Thus, by (20), 
\begin{equation} \vert \Afl\vert  \leq  C_{3} \Vert A \Vert_{L^{\infty}}. \end{equation}
The function $\ufl$ satisfies Guillemin boundary conditions along the edge
$\{x_{1}=0\}$
but with the measure $D_{0} dx_{2}$. This means that the normal component
$\Vfl^{1}$ of the vector field associated to $\ufl$ is $-D_{0}$ on the edge.
Our overall goal is, roughly speaking, to show that this is impossible if $D_{0}$ is very large. 

\


We finish this subsection with an observation which will be crucial in our
proofs below. For any two points in $Q$ the ratio of the determinant $\det
(\ufl)_{ij}$ evaluated at these two points is the same as the same as the ratio
of $\det u_{ij}$ at the corresponding points in $P$. By (21) and the observation
above we see that

\

 {\it for any two points $q,q'$on the intersection
of any   contour
$\{\frac{\partial \ufl}{\partial x_{1}}={\rm constant}\}$ with the large rectangle
$Q$  we have} 
\begin{equation} \frac{\det(\ufl)_{ij}(q)}{\det (\ufl)_{ij}(q')}< \exp(C_{1}
C_{4}). \end{equation}

 In the same vein, taking account of the rescaling, we get 
$$   (\Vfl)_{2}(x_{1}, x_{2})= \frac{ D_{0} s_{0}}{\lambda_{0}} V^{2}( 
s_{0} x_{1}, t_{0}+ \lambda_{0} x_{2}), $$
so $$ \vert \Vfl^{2} \vert \leq \frac{D_{0} s_{0}}{\lambda_{0}} C_{1}. $$
By the definition of $\lambda_{0}$ and (27) we have
$$  \frac{D_{0} s_{0}}{\lambda_{0}}= \int_{t_{0}-\lambda_{0}}^{t_{0}+\lambda_{0}}
u_{22}(\tau,0) d\tau \leq C_{4}, $$
so we have a uniform bound on $\Vfl^{2}$.

\subsection{Boundedness}

Applying the procedure above we obtain a sequence of convex  functions $\ufl^{(\alpha)}$
defined on large polygons in the half-space. In this subsection we 
show that these are bounded on compact subsets.  As before, we usually omit
the index $\alpha$ from the notation. Also, since we are only concerned with
compact subsets we can be rather vague about the precise domains of definition
of the functions.

 The first point is that, by (25), and the scalings chosen we have, for any
 fixed compact subset, a  bound
\begin{equation}  D(\ufl; p)\leq  k  \end{equation}
for points in the set, where we can suppose $k$ is as close to $1$ as we
please.

Consider the restriction
 of $\ufl$ to the $x_{1}$-axis, and set $f(x)= \ufl(x,0)$. The maximising
 condition translates into the condition that
 $$   x f'(x) - f(x) +1 \leq k x. $$
  By construction, $f(1)=f'(1)=0$ and we can integrate
 the differential inequality to obtain
 $$  f(x) \leq k x \log x .$$ 
 This discussion is valid for $x$ less than  $ s_{\delta}/s_{0}$, which we
 know is
 large.
  The normalisation conditions $\ufl(0,0)=1$
  and (27) easily imply that $\ufl(0,x_{2}) \leq 10$ (say), for $-1\leq x_{2} \leq 1$. Now
  the convexity of $\ufl$ yields a upper bound
  $$   \ufl(x_{1}, x_{2})\leq U(x_{1}, x_{2}) , $$
  say, for a suitable fixed function $U$ and all points $(x_{1}, x_{2})$
  in the triangle with vertices $(0,1), (0,-1), ( s_{\delta}/s_{0},0 )$. When
  we take our sequence $\ufl^{(\alpha)}$ the apex $s_{\delta}/s_{0}$ of this
  triangle tends to infinity.

  \

It is less easy to obtain bounds on $\ufl$ outside the triangle above, and
for this we use the bound  on the vector field $\Vfl$. As above, it suffices
to bound the function $ \ufl(0, x_{2})$. We prove
\begin{prop}
There is an $\eta>0$ such that if for some
$\sigma>0$ we have a uniform bound $ \ufl^{(\alpha)}(0, \pm \sigma )\leq U_{\sigma}$
then  there is a uniform bound $\ufl^{\alpha)}(0, x_{2})\leq
U_{\sigma+\eta}$ for all $\vert x_{2}\vert \leq \sigma+\eta$.
\end{prop}

Since $\eta$ is fixed we can use this to bound the function $\ufl$ on any
compact subset of the half-plane.

\

Write $\xi_{i}$ for the partial derivatives $\frac{\partial \ufl}{\partial
x_{i}}$. The first step in the proof of Proposition 4 is an elementary Lemma.
We consider a disc $\Delta$ of small radius $r$ centred on the point $(\frac{1}{4},
0)$.

\begin{lem}
 We can fix small $r$ and a $Z>0$ such that on the disc $\Delta$ we have
 $$   \xi_{1}<-\frac{1}{2}\ \ , \ \ \ -Z<\xi_{2}< Z, $$
 for large enough $\alpha$.
 \end{lem}
  We can choose the disc to lie well inside the triangle on which we have
  bounds on $\ufl$, and this gives bounds on the derivatives $\xi_{i}$, by
  convexity. Thus the only thing is to arrange that $\xi_{1}<-\frac{1}{2}$
  on $\Delta$. For this we begin by choosing $r$ so small that
    such that if $\vert x_{2} \vert < r$ we have
$$  \ufl(0,x_{2}) \geq .99\ \ , \ \ufl(1,x_{2})\leq .01. $$
This is clearly possible.
Consider a point $(a_{1}, a_{2})$ with $\vert a_{2}\vert <r $ and $a_{1}<1$
at which
 $\xi_{1}=-1/2$. We show that $a_{1}$ cannot be too small. By (25) we have
$$   \ufl(0,a_{2})\leq \ufl(a_{1}, a_{2})+ a_{1}/2 + k a_{1}, $$
where $k$ can be made as close to $1$ as we like.
So $$   \ufl(a_{1}, a_{2})\geq .99- (k+.5) a_{1}. $$
On the other hand, convexity implies that
$$  \ufl(1,a_{2})  \geq \ufl(a_{1}, a_{2}) - (1-a_{1})/2, $$
so $$
\ufl(a_{1}, a_{2}) \leq .5 +.01 - a_{1}/2. $$
We deduce from these inequalities that $a_{1}\geq .48 k^{-1}$ so, taking
$k$ close to $1$ and $r$ small, such a point cannot lie in $\Delta$. This establishes
the Lemma.

\

We now fix the disc $\Delta$ and the number $Z$, as above. Let $\tau_{t}$ be the
flow by translations on $\bR^{2}$
$$  \tau_{t}(\xi_{1}, \xi_{2})= (\xi_{1}, \xi_{2}+t),$$
and let $\Psi_{t}$ be the corresponding flow on the polygon $\Pfl$. That
is 
$$\Psi_{t} = (D\ufl)^{-1} \circ \tau_{t} \circ (D\ufl). $$

We consider the images $\Psi_{t}(\Delta)$, for parameters $t\in [0,T]$. Our general principle tells us that,
if these are all contained in the large rectangle $Q$, then the area of each
such image is at least a fixed multiple of the  area of $\Delta$. 
(Using the fact that the Jacobian of $D\ufl$ is $\det(\ufl)_{ij}$.)

Now suppose we have a $\sigma>r$ for which we have  obtained an {\it a priori bound}  $\ufl(0,a_{2})\leq U$, for all $\vert a_{2}\vert \leq \sigma$.
Consider a point $(x_{1}, x_{2})$ where $\xi_{1}\leq -1/2$. We first treat
the case when  $\vert x_{2}
\vert \leq \sigma$. Then we have, by convexity and the positivity of $\ufl$,
\begin{equation}   U\geq \ufl(0,x_{2}) \geq \ufl(x_{1}, x_{2}) + x_{1}/2\geq x_{1}/2.\end{equation}
Second we treat the case when $x_{2}=\sigma+t$ where $0\leq t\leq \eta$, and  $\eta$
will be chosen shortly.   
By convexity we have
\begin{equation} \ufl(x_{1}, \sigma +t) - (x_{1} \xi_{1}+ t \xi_{2}) \leq \ufl(0,\sigma)\leq U .\end{equation}
By (25) we have
$$   \ufl(0,\sigma+t)\leq u(x_{1}, \sigma+t) - \xi_{1} x_{1} + k x_{1}. $$
The two together give
\begin{equation}   \ufl(0,\sigma+t)\leq U + t \xi_{2} + k x_{1}. \end{equation}
On the other hand, the first inequality gives
$$    \ufl(x_{1}, \sigma+t) \leq U+(x_{1} \xi_{1} + t \xi_{2}),
$$ and since $\ufl(x_{1}, \sigma+t) \geq 0$ we have
\begin{equation}    x_{1} \xi_{1} \geq -(U+ t \xi_{2}). \end{equation}
Since $\xi_{1} \leq -1/2$ we obtain \begin{equation}x_{1} \leq 2(U+ t \xi_{2}).\end{equation}  Substitute back into (33) to get

\begin{equation}  \ufl(0,\sigma+t) \leq (1+2 k)( U+ t \xi_{2}). \end{equation}

Now let $\rho$ be the minimum value of $\xi_{2}$ on the set where $x_{2}=\sigma+\eta$
and $ \xi_{1}\leq -1/2$. By the above we have
$$  \ufl(0,\sigma +\eta) \leq (1+2 k) (U +\eta \rho), $$
so we want to show that $\rho$ is not too large. Let $S$ be the union of
the rectangle $\{ 0\leq x_{1} \leq 2 U, -r\leq x_{2}\leq \sigma \}$ and the
rectangle $\{0\leq x_{1} \leq (1+2 k )(U+\eta \rho) , \sigma \leq x_{2} \leq \sigma+\eta\}$.
We assume  the hypothesis that $S$ lies in the large rectangle $Q$.
By the inequalities (30),(35) above; for any $t$ with $0\leq t\leq \rho-Z$ the image $\Psi_{t}(\Delta)$ under the flow
lies in
$S$. If $\vert t-t'\vert \geq 2Z$ the images $\Psi_{t}(\Delta), \Psi_{t'}(\Delta)$ are
disjoint. Since the area of each of these is at least a fixed multiple of
the area of $\Delta$ we deduce that
$$    \rho\leq \kappa \Area(S), $$
for some $\kappa$ which does not depend on $U$. (This constant $\kappa$ depends
only on $r,Z, C_{1},C_{4}$.) Writing down the area of $S$, we get
  \begin{equation}\rho \leq \kappa \left( 2(r+\sigma)U+\eta(1+2k)(U+\eta
\rho)\right).\end{equation}
Choose $\eta$ so that $(1+2k) \kappa \eta^{2}<1/2$, say. Then we can rearrange
(37) to get an upper bound on $\rho$, which then gives by (36) a bound on $U(\sigma+\eta)$.
(Of course, we use a symmetrical argument for negative values of $x_{2}$.)
The final detail to add in the argument concerns our hypothesis that $S$
lies in the large rectangle. But, assuming this hypothesis, the set actually
lies in a smaller rectangle (because of the bound on $\rho$) so we can easily
establish the truth of the hypothesis, for large $\alpha$, by a continuity
argument.

\subsection{The blow-up limit}

Let us review the argument of this section thus far.
 We suppose we  have
a sequence $u^{(\alpha)}$ of solutions corresponding to data sets $P^{(\alpha)},
A^{\alpha)}, \sigma^{(\alpha)}$, satisfying the conditions of Theorem 3 and with $\mu_{\alpha}$
tending to infinity. We want to obtain a contradiction. We rescale
to get convex functions $\ufl^{(\alpha)}$ defined on a sequence of domains
$\Pfl^{(\alpha)}$ which exhaust the half-plane $\{x_{1}>0\}$ and we have
shown that these are
bounded on compact subsets of the closed half-plane. 
After taking a subsequence, we can suppose that the $\ufl^{(\alpha)}$ converge
uniformly on compact subsets of the {\it open} half-plane to a limit $\ufl^{(\infty)}$,
which is a continuous, weakly convex, function on the half-plane. The main
idea of our proof is to obtain a contradiction by an analysis of this limit.
First, in the next Proposition, we show that it cannot be strictly convex
anywhere. The proof we give uses the assumed bound on the associated vector fields.
The author knows of other arguments, for this step, which avoid that assumption,
but which are longer.

\

  Recall that a convex function $v$ is called strictly convex at a point $p$ if
there is an affine linear support function $\pi$ such that $v-\pi$
has a unique minimum at $p$.

\begin{prop}
The limit $\ufl^{(\infty)}$ is not strictly convex at any point of the half-plane.
\end{prop}

Suppose the contrary, so there is a small disc $D$ in the half-plane centred
at a point $p$ and an affine-linear function $\pi$ such that $\ufl^{(\infty)}-\pi$
vanishes at $p$ but is strictly positive on the boundary of $D$.  The smooth functions
$\ufl^{(\alpha)}$ satisfy elliptic equations
$$  \left(\ufl^{(\alpha)}\right)^{ij}_{ij}= - A^{(\alpha)} $$
with a fixed bound on $\Vert A^{(\alpha)}\Vert_{L^{\infty}}$. If we set
  Also
$\ufl^{(\alpha)}(p)-\lambda(p)\rightarrow 0$ and $\ufl^{(\alpha)}-\lambda\geq \delta>0$ say on $\partial
D$, once $\alpha$ is sufficiently large. Now  these facts give complete control
of the functions $\ufl^{(\alpha)}$ in the interior of $D$. We can apply Theorem
5 from \cite{kn:D3}
to obtain  upper and lower bounds on the Jacobians $\det((\ufl^{(\alpha)})_{ij})$ over the interior
of $D$ and arguing as in \cite{kn:D2}, using the theory of Cafarelli and
Gutierrez, bootstrap to control all higher derivatives.  In particular the vector fields $\Vfl^{(\alpha})$ defined by
$$  ( (V^{(\alpha)})^{i}= -\left(\ufl^{(\alpha)}\right)^{ij}_{j}, $$
are uniformly bounded on a small neighbourhood of $p$. 

To obtain a contradiction, suppose $p=(p_{1}, p_{2})$ and consider a rectangle
$$ S=\{ (x_{1}, x_{2}) : \vert x_{2} - p_{2}\vert \leq \eta, 0<x_{1}<p_{1}
\}, $$
with $\eta$ small. Since the divergence of $\Vfl^{(\alpha)}$ is bounded the
total
flux of $\Vfl^{(\alpha)}$ through the boundary of $S$ is small. The flux through
the two  edges where $x_{2}= p_{2}\pm \eta$  is bounded, by our bound on the $x_{2}$ component
of $v^{(\alpha)}$ and the flux through the edge where $x_{1}=p_{1}$ is bounded
by the argument of the previous paragraph. But the boundary conditions, after
rescaling, imply that the $x_{1}$ component of $\Vfl{(\alpha)}$ along the remaining
edge, in the boundary of the half-plane, is $D_{\alpha}$, so the flux through
this edge is $2\eta D_{\alpha}$ which tends to infinity by hypothesis.

 \

 \

 For $a\in \bR$ let $\Gamma_{a}$ denote the set of points $(x_{1}, x_{2})$
 where $x_{1}>0$ and $x_{1}= 1+ a x_{2}$. This is either a half-line or,
 in the case when $a=0$, a line.
 \begin{cor}
 The limit $\ufl^{(\infty)}$ vanishes on a set $\Gamma_{a}$ for some $a\neq
 0$.
 \end{cor}
 
 To see this let $Z$ be the zero set of $\ufl^{(\infty)}$ in the open upper
 half-plane. Recall that $\ufl^{(\alpha)}$ is normalised to achieve its minimum
 at the point $p_{0}=(1,0)$. Thus $Z$ is a convex set containing $p_{0}$.
 Proposition 5 implies that $Z$ has no extreme points and it follows immediately
 that there must be a line through $p_{0}$ whose intersection with the upper-half
 plane is contained in $Z$. We know that $\ufl^{(\alpha)}(t,0)$ is bounded
 below by $t\log t-t+1$ for $0<t<1$ and it follows that $Z$ cannot contain
 the line segment $\{x_{2}=0, x_{1}>0\}$. Thus $Z$ contains $\Gamma_{a}$
 for some $a$ and it only remains to rule out the possibility that $a=0$.
 To do this, recall that we chose our normalisation so that the $x_{2}$
 derivative of $\ufl^{(\alpha})$ differs by $1$ at the two points $(0,\pm1)$.
 This obviously implies that we can find some fixed $c$ such that $\ufl^{(\alpha)}(0,c)>2$
 say. Suppose that $Z$ contains $\Gamma_{0}$, so $\ufl^{(\alpha)}(1,
 c)$ tends to zero as $\alpha$ tends to infinity. It follows that there must be some sequence $b_{\alpha}$ tending to $1$ such that the $x_{1}$ derivatives
of $\ufl^{(\alpha)}$ evaluated at $(c,b_{\alpha})$ converge to $0$. But then
these points contradict  (30), once $k<2$.

\subsection{The final contradiction}
We again pause for discussion.
 Changing our coordinates slightly, we may without real loss
of generality suppose that $a=1$ and $\ufl^{(\infty)}$ vanishes on the ray
  $\{x_{1}=x_{2}+1, x_{1}>0\}$. The essential case to have in mind is when
  $ \ufl^{(\infty)}= \max (x_{2}-x_{1}+ 1, 0)$ so let us momentarily assume
  that we have this case.  It is tempting to try to argue
  as follows. For any large $C$ and large enough $\alpha$ we can find a point
  $p'$ near to $(C+1, C)$ such that $D(p', \ufl^{(\alpha)}$ is very close to
  $1$. In other words, transferring back to the original functions $u^{(\alpha)}$
  there are points much further from the edge than the \lq\lq worst point'' but which are \lq\lq almost as
  bad''. So this strongly suggests that if $D$ becomes large close to the
  boundary of $P$ it must also become large in the interior, which is ruled
  out by our hypotheses. Indeed if we were to drop the hypothesis on the
  integral over the boundary then we would see exactly this phenomenon, as
  we discuss further in Section 6. However, while it is suggestive, it seems
  hard to turn this line of argument into an actual proof. The proof
  we give below is rather different and hinges on 
our general principle  that $J=\det (\ufl^{\alpha})_{ij}$
changes by a bounded factor on the parts of the contours $\{\xi_{1}={\rm
constant}\}$
in the large rectangle $Q$.  We will show that this leads to a contradiction. The argument is similar to that used to prove Proposition 4 above.

\begin{lem}
The boundary values $\ufl^{(\alpha)}(0,x_{2})$ converge to $x_{2}-1$ uniformly
for $x_{2}$ in any closed interval $[-1,R]$.
\end{lem}
For given $x_{2}>-1$  we can find a sequence $x_{1}^{(\alpha)}$ converging
to $1+x_{2}$ such that the values $\ufl^{(\alpha)}$ and the partial derivatives $\frac{\partial \ufl^{(\alpha)}}{\partial
x_{1}}$ evaluated at $(x_{1}^{(\alpha)}, x_{2})$ converge to $0$ as $\alpha\rightarrow
\infty$. Then we
obtain from (30) that
$$  \ufl(0,x_{2})\leq k_{\alpha} (x_{2}+1)+\epsilon_{\alpha}$$
where $\epsilon_{\alpha}\rightarrow 0, k_{\alpha}\rightarrow 1$. By construction,
$\ufl^{(\alpha)}(0,0)=1$ and it follows from convexity that $\ufl(0,x_{2})$
tends to $x_{2}+1$, uniformly for $x_{2}$ in any compact subset of $(-1,\infty)$. However
the functions are bounded on a neighbourhood of the point $(0,-1)$ and it
follows again from convexity that the convergence is uniform up to $x_{2}=-1$.

Given a small number $r$ consider the  region 
$$\Omega =\{ (x_{1}, x_{2}): x_{1}>0, r<\sqrt{ \vert x_{1} \vert^{2}+ \vert
x_{2}+1 \vert^{2}} <r^{-1}\}. $$ 
We consider first the points in $\Omega$ where the partial derivative $\xi_{1}$ of
$\ufl^{(\alpha)}$ 
is  $-1/2$ (say) and $\xi_{2}$ is  $1/10$ (say). To simplify the exposition imagine first 
that $\ufl=\ufl^{\alpha}$ {\it vanishes} on the intersection of $\Omega$ with the ray. Thus $\xi_{1}, \xi_{2}$ also vanish on this set.
Suppose  $\xi_{1}=-1/2$ and $x_{2}>-1$. Then we must have $x_{1}< x_{2}+1$ and 
$$  0= \ufl(1+x_{2}, x_{2})\geq \ufl(x_{1}, x_{2}) -  \frac{1}{2}( 1+x_{2}-
x_{1}). $$
On the other hand
$$ (x_{2}+1) -\epsilon_{\alpha}\leq \ufl(0,x_{2})\leq \ufl(x_{1}, x_{2})
+ (k_{\alpha}+\frac{1}{2}) x_{1}, $$
where $\epsilon_{\alpha}\rightarrow 0$.  These imply that
$$   \frac{1}{2} (x_{2}+1) - k x_{1} \leq \epsilon_{\alpha}. $$
 On the other hand,
just from the fact that $\ufl(x_{1}, x_{2})\geq 0$  we have
$$ 0\leq \ufl(0,-1)+ (x_{1} \xi_{1} + (x_{2}+1) \xi_{2}). $$
So if $\xi_{1}=-1/2, \xi_{2}=1/10$ we have
$$    x_{1} \leq 2 \ufl(0,-1) + (x_{2}+1)/5. $$
If $\ufl(0,-1)$ and $ \epsilon_{\alpha}$ are sufficiently small  then these inequalities have no common solution in $\Omega$. It is clear from
a continuity argument then that the point where $\xi_{1}=-1/2, \xi_{2}=1/10$
must lie in the small half-disc $D$ of radius $r$ about the origin. 

Now obviously the same argument applies to values of $\xi_{1}, \xi_{2}$ close to $-1/2, 1/10$.
Further, it is easy to extend the argument to the case when $\ufl$ is $C^{0}$
close to a function vanishing along the ray, over the fixed annulus. So we
conclude that  there is a small rectangle $R\subset \bR^{2}$ of the form
 $$R=\{ (\zeta_{1}, \zeta_{2}): \vert \zeta_{1}+1/2 \vert<\eta, \vert \zeta_{2}-1/10\vert <\eta\} $$
with the following property. 
 For any given $r$ and all large enough $\alpha$   all points $(x_{1},
 x_{2}) $ for which
 which  $(\xi_{1}(x_{1}, x_{2}),
 \xi_{2}(x_{1}, x_{2}))\in R $ are contained in  $D$. 
 
 Now for fixed large $\alpha$ take a point $(\zeta_{1}, \zeta_{2})$ in $R$ and
 consider the contour $\xi_{1}(x_{1}, x_{2})=\zeta_{1}$. This contour meets the line $x_{2}=0$
 at some point $(a_{1}, 0)$. As in the proof of Proposition 4 we have $a_{1}>c$ for some fixed $c>0$.  Then using the estimate in Theorem 5 of \cite{kn:D3} we have
 an upper bound on the determinant function at this point. It is obvious
 that when $x_{1}>-1$ the contour cannot move out of the large rectangle
 $Q$. We conclude from
 our principle that
 the determinant is bounded at all points whose derivative lies in $R$, say
 $J\leq C$. But the inverse of the derivative maps $R$ into $D$ so
 $$  \int_{R} J^{-1} d\zeta_{1} d\zeta_{2} \leq \Area(D)= \pi r^{2}/2. $$
 
 Thus $4 C^{-1} \eta^{2} \leq \pi r^{2}$. But since $r$ can be made arbitrarily
 small, with $\eta$ fixed, this gives our contradiction.
 
 \input extremal2b.tex

%% file: extremal2b.tex
\section{$C^{\infty}$ limits away from the vertices}
In the previous section we obtained a uniform bound on the quantity $D(p;u)$
along the interior of each edge. We now use this to get complete control
of the solution away from the vertices. This is straightforward, given the
results from \cite{kn:D3}, if we have
a lower bound on the Riemannian distance function determined by the solutions,
and we explain this argument in (4.1). The main work of the section goes
in to establishing this lower bound. For this we derive various estimates
on the solution, near to an edge, and particularly on 
$\det u_{ij}$. These estimates may have independent interest.  

\subsection{The proof, assuming a lower bound on the Riemannian distance}

We begin with the relation between the quantity $D(p;u)$ and the \lq\lq M-condition''.
Recall that in \cite{kn:D3} we said that $u$ satisfies an $M$ condition if $V(p,q)\leq
M$
for any pair $p,q$ of points in $P$ such that the line segment
$ \{t p + (1-t) q: -1\leq t \leq 2 \}$
lies in $P$. Here $V(p,q)$ is the variation of the derivative of $u$ in
the direction of the unit vector $\nu= (p-q)/\vert p-q\vert$ between 
the two points. For brevity we will call such pairs $p,q$ \lq\lq admissible
pairs''.

\begin{prop}
Suppose $u$ is a normalised function on the polygon $P$ and we have
\begin{itemize}
\item A bound on the integral of $u$ over $\partial P$;
\item For each $\delta>0$ a  bound on $D(p;u)$ for points $p$ whose Euclidean
distance to all vertices of $P$ exceeds $\delta$.
\end{itemize}
Then for any $\delta'>0$ the variation $V(p,q)$ is bounded for all admissible pairs $p,q$ where the Euclidean distance from  $p$ to  the vertices exceeds
$\delta'$.
\end{prop}

This is very elementary, so we will use rather informal language. The first
hypothesis controls $V(p,q)$ when $p$ is not close to the boundary,  so the relevant case is when $p$ is close to a unique edge. We suppose,
as in the previous section, that this edge is a segment
$a\leq x_{2}\leq b$ of the $x_{2}$-axis and $p=(p_{1}, p_{2})$ with $a+\delta<p_{2}<b-\delta$.
Let $u_{p}$ be the function obtained from $u$ by normalising at $p$. As in
Proposition 4, the
bound on $D(p')$, for points $p'$ on the segment $\{p'_{2}=p_{2}\}$, gives a bound
\begin{equation}   u_{p}(p'_{1}, p_{2})\leq  C p_{1}, \end{equation}
for $0\leq p'_{1}\leq 3 p_{1}$, say. With this point $p$, the points $q$
we need to consider in the definition of the $M$-condition range over some
quadrilateral $Q$. Two of whose edges are segments in the lines
 $\{x_{1}=p_{1}/2\},\{x_{1}=2p_{1}\}$
and the other two are determined by the other edges of $P$. But these other two
edges are a definite distance from the rest of the boundary of $P$. We can
choose a slightly larger quadrilateral $Q^{+}$, two of whose edges are segments
in the lines $\{ x_{1}=0\}, \{x_{1}= 2p_{1}\}$ and whose other two edges
are again a definite distance from the rest of the boundary of $P$. For each
unit vector $\nu$ there are unique $h,h^{+}$ such that $q=p+h\nu$ lies in the
boundary of $Q$ and $q^{+}=p+h^{+} \nu$ lies in the boundary of $Q^{+}$. We can
suppose that $h\leq (1-\epsilon) h^{+}$ for some fixed $\epsilon>0$. Then
if $q=p+h\nu$ we have
$$   V(p,q)\leq \frac{u_{p}(q^{+})}{\epsilon h^{+}}. $$
Write the $x_{2}$ coordinate of $q^{+}$ as $p_{2}+t$. When $t=0$ then (38) states
that
$u_{p}(q^{+})\leq C p_{1}$. When $q^{+}$ lies on one of the other two edges
of $Q^{+}$ (not parallel to the $x_{2}$-axis) we have a bound $u_{p}(q^{+})\leq
C$, since these edges are a definite distance from the other edges of $P$.
Convexity of $u_{p}$ yields an inequality of the form
$ u_{p}(q^{+}) \leq C (p_{1} + \vert t\vert)$.  On the other hand we have
$h^{+}\geq C \sqrt{ p_{1}^{2}+t^{2}}$ so
$$   \frac{u_{p}(q^{+})}{ h^{+}} \leq C \frac{p_{1}+\vert t\vert}{\sqrt{p_{1}^{2}+t^{2}}},
$$ which is bounded. This completes the proof.

\

\

Let $\Omega\subset \oP$ be the set obtained by deleting fixed small Euclidean
discs
about the vertices and let $\partial_{*}\Omega$ be that part of the boundary
of $\Omega$ which is {\it not} contained in the boundary of $P$. An admissible convex function $u$ on $P$ defines a Riemannian
metric $u_{ij}$on $\oP$, regarded as a $2$-manifold with corners. For $p,q\in
\Omega$ we write
$\distu(p,q)$ for the Riemannian distance defined by this metric, and
$$  \distu(p,\partial_{*}\Omega) = \inf_{q\in \partial_{*}\Omega} \distu(p,q).
$$

 Locally, we may also associate a $4$-dimensional Riemannian
manifold to this data and we let $\vert F\vert^{2}$ be the square of the
Riemannian norm of the curvature tensor as in \cite{kn:D2}, \cite{kn:D3}. Now consider
a sequence  $u^{(\alpha)}$ as in Theorem 3. We claim  that
\begin{prop}
There is a fixed bound $\vert F^{\alpha}\vert^{2} {\bf dist}(\ , \partial \Omega)^{2}\leq C$
for all $\alpha$.
\end{prop}
Of course here, strictly speaking we have a sequence of polygons $P^{(\alpha)}$
so we need to fix a sequence of domains $\Omega^{(\alpha)}$, but the meaning
should be clear.

The proof of Proposition 7 is a straightforward modification of  the arguments of \cite{kn:D3},
which we only outline. We proceed by contradiction
and suppose there is a sequence of points $p_{\alpha}$ for which $K_{\alpha}=\vert F\vert^{2}
{\bf dist}(\ ,\partial_{*}\Omega)^{2}$ tends to infinity. Then we rescale the
metric so that after rescaling the  curvature has norm $1$ at the chosen
points. After this rescaling the distance to the boundary $\partial_{*}\Omega$ is $K_{\alpha}$
which becomes large by hypothesis, and the curvature is bounded on balls of
a fixed size about the chosen points. This means that we can take the blow-up limit just as in \cite{kn:D3} and the extra boundary \lq\lq disappears'' in the limit.
Then the analysis of the blow-up limits in \cite{kn:D3} gives the desired
 contradiction. We have to use the $M$-condition a number of times in these
 arguments, but only at points in $\Omega$, and we have this  by Proposition
 6.

 Now fix a subset $\Omega_{0} \subset \Omega$, for example removing larger Euclidean discs
 about the vertices. We will show
 \begin{prop} There is an $\eta>0$ such that, for all $\alpha$
 $$  \distual( \partial_{*}\Omega_{0}, \partial_{*}\Omega)\geq \eta. $$
 \end{prop}
 Assuming this,  Proposition 7 gives  an upper bound on the size of the curvature tensor over $\Omega_{0}$
 and the arguments of \cite{kn:D3} apply without change to give $C^{\infty}$
 convergence of the $u^{(\alpha)}$, in the same sense as in
 \cite{kn:D3}. Since we can make $\Omega, \Omega_{0}$ as large as we please,
 we conclude that the $u_{\alpha}$ converge in $C^{\infty}$ on compact subsets
 of $\oP$ minus the vertices. The proof of Proposition 8 takes up the remainder
 of this section.
 
 \subsection{Lower bound on Riemannian distance: strategy}
 
 By the results of \cite{kn:D2} we know that over any compact subset of the
 open polygon $P$ the Riemannian length of paths compares uniformly with
 the Euclidean length. Also we know that, given a bound on the quantities
 $D(p)$, the Riemannian
 length of a line segments meeting an edges in an interior points is bounded
 below (\cite{kn:D3}, Lemma 2). Using these facts, it is elementary to reduce
 the proof of Proposition 8 to the following. Given any two points $q, q'$ in
 the interior of an edge, there is a lower bound on the Riemannian length
 of paths from $q$ to $q'$. We can take the edge to be a segment in the
 $x_{2}$-axis and $q=(0,\alpha), q'=(0,\beta)$, with $\alpha<\beta$. The same elementary arguments
 show that it suffices to consider paths which lie in a rectangle
 $\{(s,t): 0\leq s\leq s_{0}, \alpha\leq t\leq \beta\}$, for arbitrarily small
 $s_{0}$. 
 
 {\bf Remarks}
 \begin{enumerate}
 \item Of course we fix $s_{0}$ so that this rectangle is well away from the other
 edges of $P$.
 \item  It is not hard to avoid appealing to the results of \cite{kn:D2} here, at
 the cost of some extra arguments.
 \item The obvious path, given by the line segment in the $x_{2}$-axis,
 between these points is a geodesic and we {\it expect} that this will be
 the length minimising path. If we knew this then the proof of Proposition
 8 would be
 substantially simpler---we could avoid Proposition 10 below---but the author
 has not found a argument to establish this fact so we have to work more.
 
 \end{enumerate}
 
 \
 
 Our basic idea is to consider the function $\xi_{2}= \frac{\partial u}{\partial
 x_{2}}$ on $P$.  For a pair
 of points $(0,t_{1}), (0,t_{2})$ on the edge, with $t_{1}<t_{2}$, write $$\Delta(t_{1}, t_{2})= \xi_{2}(0,t_{2})-\xi_{2}(0,t_{1})=\int_{t_{1}}^{t_{2}}
u_{22}\ dt .$$
One step in the proof is to establish that $\Delta(\alpha,\beta)$ is not
small (Corollary 4 below). To see the relevance of this consider, for this exposition, the linear path along
the $x_{2}$-axis. The Riemannian length of this path is
$$  \int_{\alpha}^{\beta} \sqrt{u_{22}} \ dt, $$
while $$ \Delta(\alpha,\beta)= \int_{\alpha}^{\beta} u_{22} \ dt. $$
Informally, we expect that if $\Delta(\alpha,\beta)$ is not small then
$u_{22}$ should not be small at typical points and so the Riemannian length
should not be small. But of course this argument does not suffice, as it stands, because
of the square-root in the integral for the Riemannian length. Much the same
issue arose  in \cite{kn:D2}, deriving estimates in the interior of the polygon.
The analogous difficulty there was to obtain lower bounds for the Riemannian
distance to the boundary given a \lq\lq strict convexity'' condition. The
approach in \cite{kn:D2} was to replace  Riemannian balls  with  \lq\lq
sections''of the convex function, using deep results of Caffarelli. The problem
at hand is that we are working up to the boundary, where these results do
not apply.

To proceed with our outline of the strategy, consider the square of the Riemannian norm
 of its derivative $\nabla \xi_{2}$ which is
 $$ \vert \nabla \xi_{2} \vert^{2} = u^{ij} \sum \frac{\partial \xi_{2}}{\partial
 x_{i}}\frac{\partial \xi_{2}}{\partial x_{j}} = u^{ij} u_{2i} u_{2j} = u_{22}.
 $$
 So along any path in $P$ with the given end points the change in $\xi_{2}$ is bounded by $\int \sqrt{u_{22}}
 \ d\sigma$, where $d\sigma$ denotes Riemannian arc length along the path.
 Thus if we have an {\it upper} bound $u_{22}\leq C$ along the path, we have
 $$ \Delta(\alpha,\beta)\leq L \sqrt{C}, $$
 where $L$ is the length of the path; so we have a lower bound $L\geq C^{-1/2}
 \Delta(\alpha, \beta)$, as desired. This upper bound on $u_{22}$ is given in Proposition 10 below. The proof of this, and the lower bound on $\Delta(\alpha,
\beta)$ goes through estimates for the determinant $J=\det u_{ij}$. We emphasise
that in all of these arguments we make
much use of
the result of Section 4; $D(p;u)\leq D$ say, for all relevant points $p$,
and the various constants in our statements depend on $D$.

\subsection{Lower bound on Riemannian distance: detailed proofs}
 
 \begin{lem}
 Let $t_{1}<t_{2}$ be two points in the interval
 $[\alpha,\beta]$ and $\tau=(t_{1}+t_{2})/2$. 
 There are constants $c_, c'$  such that if for some $s\leq
 s_{0}$
 we have  $\Delta(t_{1}, t_{2})\leq c  \frac{s}{t_{2}-t_{1}}$
 then $J(s,\tau)\leq c' (t_{2}- t_{1}) ^{-2}$.
 \end{lem}
 The proof is sufficiently like Lemma 14 in \cite{kn:D3} that we leave this to the
 reader. (Elementary arguments give bounds on the first derivative in a suitable
 neighbourhood of $(s,\tau)$, then we apply the maximum principle result Theorem
 5
 of \cite{kn:D3}.)

 \begin{cor}
 For any $\mu>1$ there is a constant $C_{\mu}$ such that
 $$  \Delta(t_{1}, t_{2}) \geq C_{\mu} (t_{2}-t_{1})^{\mu}. $$
 In particular $\Delta(\alpha, \beta)$ is bounded below.
 \end{cor}
 
 To see this use Theorem 5 in \cite{kn:D2} which states that for any $a<1$ the function
 $J$ satisfies a lower bound $J(s,t)\geq C s^{-a}$. Then the statement follows
 immediately after re-arranging the inequalities. Note that if we could take
 $\mu=1$ we would we in a strong position---$u_{22}$ would then be bounded
 below on the edge--- but the author has not been able
 to achieve this directly.

 The next step is to   find a sharp {\it upper} bound on the function $J$.

 \begin{prop}
 There is a constant $C$ such that $J(s,t)\leq C s^{-1}$ for all
 $t\in [\alpha, \beta]$.
 \end{prop}
 We choose nested intervals $(\alpha, \beta)\subset (\alpha', \beta')\subset
 (\alpha'', \beta'')$ so that the rectangle
 $(0,s_{0}]\times [\alpha'', \beta'']$ is well away from the other edges
 of  $P$. We have an upper bound  $\Delta(\alpha'', \beta'')\leq \Delta''$
 say. We can suppose that the result of Lemma 8 applies in the larger
 interval $[\alpha'', \beta'']$. We have an 
 upper bound, $J(s_{0},t)\leq J_{0}$ say, if  $\alpha''<t<
 \beta''$. Set $\eta_{0}= J_{0}^{-1} s_{0}^{-1}$ and consider some $\eta$ with
 $0<\eta\leq\eta_{0}$. 
 
 \
 
Consider the function $F= J^{-1}$. This satisfies the equation $u^{ij} F_{ij}
= -A<0$ (see (14) in \cite{kn:D2}). Thus  the function $G=F-\eta x_{1}$ has no interior minima. We have
$G=0$ on the axis $\{x_{1}=0\}$ and $G>0$ on the parallel segment $\{x_{1}=s_{0}, \alpha''<x_{2}<\beta''\}$. Let
$Q$ be the rectangle $(0,s_{0})\times (\alpha', \beta')$ and $\Sigma$ be the
subset of $Q$ on which $G<0$.
 Suppose $\Sigma$ contains a point $p=(p_{1}, p_{2})$ with $\alpha<p_{2}<\beta$. Then the connected
 component of $\Sigma$ containing $p$ must meet the boundary of $Q$, since there
 are no interior minima and by construction this can only occur on the boundary
 components $x_{2}=\alpha', \beta'$. So there is a continuous path in $S$ from $p$
 to either the  boundary $x_{2}=\alpha'$ or to $x_{2}=\beta'$. Without loss
 of generality suppose the former. Then for each $\tau\in (\alpha', \alpha)$
 there is a point $(s, \tau)$ in $\Sigma$, {\it  i.e.} where $J(s,\tau)> \eta^{-1} s^{-1}$.
 Now let $\lambda= \sqrt{c' \eta s}/2$, with $c'$ as in Lemma 8.  We suppose
 $\eta$ is chosen so that $\sqrt{c'\eta s_{0}}/2 < \alpha'-\alpha''$ thus
 $\lambda<\alpha'-\alpha''$ and the interval $[\tau-\lambda, \tau+\lambda]$
is contained in $[\alpha'', \beta'']$. By Corollary 4 we have
\begin{equation} \int_{\tau-\lambda}^{\tau+\lambda} u_{22} dt \geq c \frac{s}{\lambda} = \left(\frac{4c}{c'\eta}\right) \lambda . \end{equation}

\

Let $f$ be the restriction of the second derivative $u_{22}$ to the interval
$[\alpha'', \beta'']$ in the edge, extended by zero to a function on
$\bR$. Thus  $\Vert f \Vert_{L^{1}}\leq \Delta''$.
Let $m_{f}$ be the {\it maximal function} of $f$;
$$  m_{f}(\sigma)= \max_{\mu>0} \frac{1}{2\mu} \int_{\sigma-\mu}^{\sigma+\mu}
f(t) \ dt. $$Thus
$$ \frac{1}{2\lambda}\int_{\tau-\lambda}^{\tau+\lambda} u_{22} dt \leq m_{f}(\tau).
$$
and (39) gives $m_{f}(\tau) \geq \frac{2c}{2c'\eta}$ for each $\tau \in [\alpha',
\alpha]$. Now the weak type bound for the maximal function tells us that there
is a constant $C$ such that for all $b$ the
measure of the  set on which $m_{f}$ exceeds $b$ is at most $C \Vert f \Vert_{L^{1}}
b^{-1}$. Thus
$$  (\alpha-\alpha')\leq \frac{C\Delta'' c'}{2c} \eta.
$$
If we choose $\eta$ sufficiently small we get a contradiction, so there can
be no such point $p$. In  other words $J(s,t)\leq \eta^{-1} s^{-1}$
for $t \in [\alpha, \beta], s\leq s_{0}$.

\

It is easy to see, from the Guillemin boundary counditions, that the limit
of $sJ(s,t)$ as $s\rightarrow 0$ is the second derivative $u_{22}$, evaluated
at the point $(0,t)$. So a Corollary of the result above is that $u_{22}$
is bounded on the interval $[\alpha, \beta]$ in the edge. This then gives
us a lower bound on the Riemannian length of this interval and, as in the third
remark at the beginning of Section 4.1, we strongly suspect that this actually
realises the minimal length. However, lacking a proof of this, we go on to
prove.
\begin{prop}
There is a constant $C$ such that $u_{22}\leq C$ at all points $(s,t)$ with
$s\leq s_{0}, \alpha\leq t\leq\beta$.
\end{prop}
This result completes the proof of Proposition 8, as we explained in 4.1.

  \
  
  The proof of Proposition 10 is roughly speaking to argue that if $u_{22}$
  is large then $J$ would violate the bound of Proposition 9.

\begin{lem}
Given $D>0$ there are positive $\kappa, \zeta_{1}, \zeta_{2}, R>1$ with the following property.
Suppose $v$ is a smooth convex function on th rectangle $\{0\leq x_{1}\leq R , -R\leq
x_{2}\leq R\}$,  whose derivative is a diffeomorphism to its
image. Write $v_{i}$ for the partial derivatives $\frac{\partial
v}{\partial x_{i}}$. Suppose that $v_{1}(x_{1}, x_{2})\rightarrow -\infty$
as $x_{1}\rightarrow 0$. Suppose that  $v$ satisfies a bound $D(p;v)\leq D$ for all points $p$.
 Suppose that $v$ is normalised at the point
$(1,0)$, that $v(1,x_{2})\leq 2D$ for $-1\leq x_{2}\leq 1$ and $v(1,1)=2D$.
Then any point $(x_{1}, x_{2})$ where $v_{1}\leq -\zeta_{1}$ and $\zeta_{2}<v_{2}<2\zeta_{2}$
has $\vert x_{1}\vert, \vert x_{2} \vert \leq \kappa$.
\end{lem}

The proof of this is similar to the arguments in Proposition 4 and Lemma
7. All the steps are
entirely elementary so we will use informal language.
Given $\zeta_{i}$, let $S$ be the set of points with $v_{1}\leq -\zeta_{1}$ and $\zeta_{2}<v_{2}<2\zeta_{2}$. We choose $\zeta_{2}\leq 1/2$ so the hypotheses
imply that for any $\zeta\in [\zeta_{2}, 2\zeta_{2}]$ there is an $x_{2}$
in $(0,1)$ such that $v_{2}=\zeta$ at the point $p_{\zeta}=(1,x_{2})$. We consider the
contour $\Gamma$ on which $v_{2}=\zeta$ and $x_{1}\leq 1$. This cuts each line $\{x_{1}={\rm constant}\}$
exactly once (so long as the intersection point does not move out to the boundary
$x_{2}=\pm R$). We have $v(x_{1}, 0) \leq D$ for $0\leq x_{1}\leq 1$.
Then convexity implies that no point $(x_{1}, x_{2})$ with $0\leq x_{1}\leq
1$ and $x_{2}$ very negative can lie in $\Gamma$.
Given  a large positive $\rho$ we consider the line through the points $(1,2)$ and
$(0,\rho)$. We have an upper bound on the value of $v$ at the intersection of this line with the $x_{1}$
axis while $v(1,2)=2D$ by hypothesis. Convexity implies that at a point
$(x_{1}, x_{2})$ on this line with $0\leq x_{1}\leq 1$ the value of $v$ must
be approximately $2 D x_{2}$. For suitable choices of the parameters we see
that the contour $\Gamma$ cannot meet this line segment. In particular, $\Gamma$
is confined to lie in a bounded region $Q=\{0\leq x_{1}\leq 1, -\rho\leq x_<{2}\leq \rho\}$ say. (We can suppose $R>\rho$ so we do not have any difficulties with
the domain of definition. ) As we move along the contour $\Gamma$, with $x_{1}$
decreasing, the derivative $v_{1}$ tends to $-\infty$ so whatever the value of
$\zeta_{1}$ the point on the contour eventually lies in $S$. On the other
hand if we choose $\zeta_{1}$ large then the bound on the  $D(p)$ implies that
$p_{\zeta}$ is {\it not} in $S$. The hypotheses imply that $S$ is connected
and it follows that $S$ is contained in the bounded set $Q$, which completes
the proof.

We now prove Proposition 10.  Given a point $p=(s,t)$, we  let $u_{p}$ be
the function obtained by normalising $u$ at $p$. By applying 
Theorem 5 in \cite{kn:D3} together with lower bound on $J$, much as in the proof
of Lemma 8, we  see that there is a small positive number $\mu$ such that either
$u_{p}(s,t+\mu)= D s$ or $u_{p}(s, t-\mu)=D s$. Without loss of generality
suppose the former, and that $\mu$ is the least possible such value. Write
$\mu=r \sqrt{s}$. Now  rescale to define

\begin{equation}   \ufl(x_{1}, x_{2}) = s^{-1} u_{p}(s x_{1},r \sqrt{s} x_{2}+t).\end{equation}
Then $\ufl$
satisfies the hypotheses on the function $v$ of Lemma 9 (and we can suppose
$R$ is as large as we please, since we are only concerned with small $s$) We write $\xifl_{1}, \xifl_{2}$ for the derivatives of $\ufl$. We see that from Lemma 9 that points with
$\xifl_{1}<-\zeta_{1}$ and $\zeta_{2}<\xifl_{2}< 2\zeta_{2}$ lie in a fixed bounded set.

Now write $\Vfl$ for the vector field associated to $\ufl$,  as in Section
4. Calculating the transformation under rescaling (40) we find that $\Vfl$ is
{\it bounded}. Let $\phi$ be the Legendre transform of $\ufl$ and consider
the rectangle
$$  Q=\{ (a_{1}, a_{2}): -(\zeta_{1}+1)\leq a_{1} \leq a_{2} , \zeta_{2}\leq
a_{2} \leq 2 \zeta_{2}\}. $$
The bound on $\Vfl$ means that the determinant of the Hessian of $\phi$
varies by a bounded factor over $Q$. Since this derivative maps $Q$ into
a bounded set we get an upper bound on this determinant at each point of
$Q$. Further, for any given $\rho$ we get an upper bound on  the Hessian
over the whole ball $\vert \ua \vert \leq \rho$. Now the choice of scaling,
and the bound on $D(p, )$,
gives bounds on the derivative of $\ufl$ over a disc of radius $1/4$, say,
centred at $(1,0)$. Since the determinat of the Hessian of $\phi$ is the
inverse of $\det (\ufl)_{ij}$, at the corresponding point, we obtain a {\it lower} bound
on $\det (\ufl)_{ij}$ over this disc. But we also have an upper bound on
this determinant, by Lemma 14 of \cite{kn:D3}. Then we deduce, just as in \cite{kn:D2}, bounds
on  all derivatives of $\ufl$ on a small neighbourhood of the point $(1,0)$
 In particular $$\vert\frac{\partial^{2} \ufl}{\partial x_{2}^{2}}\vert \leq c $$ and
$$ \det((\ufl)_{ij}) \geq c^{-1} $$
say.
Now we have the transformation relations, from (40), 
$$    (\ufl)_{22}= r^{2} u_{22}\ ,\ \det((\ufl)_{ij})= r^{2} s \det(u_{ij}). $$
Since $s \det(u_{ij})\leq C$ by Proposition 9, we deduce that
$$u_{22}\leq c^{2} C$$ as required.

\section{The vertices}

Let us again take stock of our progress.  We are considering a convergent
sequence of data sets $(P^{(\alpha)}, A^{(\alpha)}, \sigma^{(\alpha)})$ with
solutions $u^{(\alpha)}$ normalised at the centre of mass of  $P^{(\alpha)}$.
Our original hypothesis is that the integrals of $u^{(\alpha)}$ over $\partial
P^{(\alpha)}$ are bounded, and we showed in Section 2 that the $u^{(\alpha)}$
satisfy an $L^{\infty}$ bound. Then we saw in Sections 3 and 4 that
the $u^{(\alpha)}$ converge away from the vertices. Our task in this section
is to show that the solutions converge in neighbourhoods of the vertices. We can fix attention on a single vertex and we choose
coordinates so that this vertex is the origin, that  $P=P^{(\alpha)}$ is
equal to the quarter
plane $\{ x_{1}, x_{2}>0\}$ near the vertex and the measures on the two edges
$\{x_{1}=0\}, \{x_{2}=0\}$ are standard. As before we usually omit the index
$\alpha$. 
For sufficiently small positive $t$ we write
\begin{equation} E(t)= t^{-1}( u(2t,0)+ u(0,2t)- 2 u(t,t)). \end{equation}
 Our strategy
is to prove
\begin{prop}
There is a bound $E(t)\leq E_{0}$ for all $t, \alpha$.
\end{prop}
Of course, this is only  of interest  for small values of $t$.
Given this, it is not very difficult to deduce the desired convergence around
the vertex, see subsection 5.5.

Our proof of Proposition 11 is complicated, so we will first give some discussion to motivate
the constructions.  The bound is similar in character to the bound on the
quantity $D$ which we obtained in Section 3, and some of the same difficulties
emerge in the proof. For each $\alpha$ choose a value $t_{0}$ which maximises the function
$E$ and set $\Emax= E(t_{0})$. Define a function $\ufl=\ufl^{(\alpha)}$ by
$$  \ufl(x_{1}, x_{2})= \Emax^{-1} t_{0}^{-1} \left( u(t_{0}x_{1}, t_{0} x_{2}) +\pi(x_{1},
x_{2}) \right) $$
where $\pi$ is the affine-linear function chosen so that $\ufl$ is normalised
at the point $(1,1)$. 
We suppose that, in the sequence $(\alpha)$, the maxima $\Emax=\Emax^{(\alpha)}$ tend to infinity and seek a contradiction. It is not hard to show that  the $\ufl$ converge to a convex function $\ufl^{(\infty)}$ but the main difficulty
is to rule out the possibility that
$$   \ufl^{(\infty)}(x_{1}, x_{2})= \frac{1}{2} \vert x_{1}- x_{2} \vert. $$
Compare with the discussion in (3.4) above, for the quantity $D$. To get around
this we consider also the determinant function $J=\det u_{ij}$ and make various
arguments with this. A crucial point is that, using the $L^{\infty}$ bound
from Section 2, we obtain sharp upper and lower bounds on $J$ in terms of
the Legendre transform coordinates $\xi_{i}$ (Proposition 12 below). Then
we consider a \lq\lq perturbation'' of the function $E$ and maximise this
to obtain, ultimately, the desired contradiction. (In fact we do not explicitly
pass to  the limit $\ufl^{(\infty)}$ in our actual proof, making all our
arguments with the smooth functions $u^{(\alpha)}$, but the reader may find
it helpful to have this in mind when following the arguments.)

\subsection{Volume bound}

We continue with the same notation reviewed above, focussing on a vertex
$(0,0)$ and, given $u$, we set $\xi_{i}= \frac{\partial u}{\partial x_{i}}$.
We write $J=\det u_{ij}$. Notice that for the flat model we have
\begin{equation}J= (x_{1} x_{2})^{-1}= e^{\xi_{1} + \xi_{2}} \end{equation}

\

\begin{prop}
There is a constant $B$ such that
$$   B^{-1} e^{\xi_{1}+ \xi_{2}} \leq J \leq B e^{\xi_{1}+\xi_{2}}
$$ in a fixed neighbourhood of the vertex.
\end{prop}

Fix some standard reference sympletic potential function $u_{0}$ (so really
we have a convergent sequence $u_{0}^{(\alpha)}$). Let $\phi, \phi_{0}$ be the Legendre
transforms of $u,u_{0}$ respectively. We have an elementary identity
$$   \Vert \phi-\phi_{0} \Vert_{L^{\infty}} = \Vert u - u_{0}\Vert_{L^{\infty}}.
$$
Clearly the $u_{0}$ are bounded and so by Theorem 2 the difference $\phi-\phi_{0}$
is bounded.  Now take complex coordinates $z_{1}, z_{2}$ and set $\xi_{i}=
\log \vert z_{i}\vert$, so we regard $\phi, \phi_{0}$ as functions of the
$z_{i}$. Fix a neighbourhood $N$ of the vertex in $\oP$. Under the Legendre
transform, this corresponds
to some neighbourhood $U$ of the origin in $\bC^{2}$. The results of the
previous section give upper and lower bounds on the difference $\log J- (\xi_{1}+\xi_{2})$  over the boundary of $U$. Since the origin is a vertex of the polygon, these functions extend
to smooth functions on $\bC^{2}$. The function $\phi_{0}$ satisfies some
fixed bound on the unit ball $B^{4}\subset \bC^{2}$ so, by the above, $\phi$
does also. 

The results of the previous sections give us $C^{\infty}$ bounds on  $\phi$
over compact subsets of the punctured ball $B^{4}\setminus \{0\}$. 
Let $V$ be the  volume element of the metric  in these complex co-ordinates,
that is $V=  \det (\frac{\partial^{2}\phi}{\partial z_{i}\partial \overline{z}_{j}})$.
So we have an upper and lower bounds on $V$ away from the origin in $B^{4}$.
The prescribed scalar curvature equation is
$$ \Delta \log V = A, $$
where $A$ is thought of as a function on $\bC^{2}$ via the Legendre transform
and $\Delta$ is the usual Laplace operator of the Kahler metric.
Thus $ \vert \Delta \log V \vert \leq C $ say. Since $\Delta \phi=2$ we have
$$ \Delta (\log V+ \frac{C}{2} \phi) \geq 0\ \ ,\ \  \Delta (\log V-\frac{C}{2} \phi )\leq 0.$$
Thus, by the maximum principle and our bound on $\phi$, the function $\log
V$ over
the entire ball is controlled by its values on the boundary, so we have
upper and lower bounds on $V$ over $B^{4}$. Now the chain rule gives
$$  \det u_{ij} = V^{-1} \exp(\xi_{1}+ \xi_{2}) $$
and our result follows.

\

Next we have a simple  lower bound on the determinant $\det u_{ij}$.
\begin{lem}
There is a constant $c>0$, depending only on $B$ above , such
that $ J\geq c (x_{1}+ x_{2})^{-2}$.
\end{lem}

To see we argue in the same manner as in Lemma 3. We consider a point $\up=(p_{1},
p_{2}) $ in the quadrant $\{x_{1}, x_{2}>0\}$ and let $u_{p}$ be the function
obtained from $u$ by normalising at $p$.
 Let $Q$ be the square consisting of points $(\xi_{1}, \xi_{2})$
with $\vert \xi_{i}+1 \vert \leq 1/10$ (say) and let $S$ be the set of points $(x_{1},
x_{2}) $ at which the derivative of $u$ lies in $Q$. The previous result implies that over $S$
$J$ differs by a bounded factor from $J(p)$. So we have
\begin{equation}  \Area(S)= \int_{Q} J^{-1} d\xi_{1} d\xi_{2} \geq c J(p)^{-1}. \end{equation}
 Let $\uy$ be a point of $S$ and $\pi$ be the affine-linear function defining the supporting hyperplane at $\uy$.
The zero set of $\pi$ is a line which separates $\uy$ and $p$ and it follows
from this that $\uy$ lies in the triangle
with vertices $(0,0), (p_{1}+ \frac{11}{9}
p_{2},0), (0, p_{2} +\frac{11}{9} p_{1}$. So   the area of $S$ is not more
than $(p_{1} +p_{2})^{2}$. Rearranging (43) then gives the result.

Notice that, comparing with (42), the bound in Lemma 10 is in a sense sharp when $p_{1}, p_{2}$ are approximately
equal.

\subsection{Proof on the diagonal}
Recall the definition of $E(t)$ in (41). In this subsection, and the next two, 
 we find an {\it a priori} upper  bound on $E(t)$.
\begin{prop}
There is a bound $E(t)\leq E_{0}$ for all $t, \alpha$.
\end{prop}
Notice that $E(t)$ is not changed if we add an affine-linear function to
$u$ and  that $E(t)$ is preserved by the rescaling
$$   \tu(x_{1}, x_{2}) = \lambda^{-1} u(\lambda x_{1}, \lambda x_{2}). $$
Under this rescaling the function $A$ transforms to $\lambda A$.
Making this rescaling, with small $\lambda$, and changing notation in the
obvious way, we can suppose that $u$ is defined on a large region in the
quarter-plane $\{x_{i}>0\}$. It seems simplest to take this rescaling
as understood, without bringing in explicit notation. By the scaling behaviour,
we can suppose that $\Vert A\Vert_{L^{\infty}}$ is as small as we please:
let us suppose it is less than $1$.  It will often be convenient to work in the co-ordinates
$$t=\frac{1}{2} (x_{1}+x_{2}), s= \frac{1}{2} (x_{1}-x_{2}). $$
 Recall that we set
$\Emax= {\rm max}_{t>0} E(t)$. We also write $J(t)$ for the determinant
$\det( u_{ij})$ evaluated at $(t,t)$ and we write $u(t)$ for the function
of one variable $u(t,t)$. For integers $n\geq 0$ let 
$$  \delta_{n}= u'(2^{-n+1}) - u'(2^{-n}). $$
\begin{prop}
There is a  constant $c$ such that then
$\delta_{n}\leq 2 \log \Emax  +c $ for all $n$,
\end{prop}

To prove this we observe that $E(1)$ controls the variation in the partial
derivative $\frac{\partial u}{\partial s}$ over an interval in the line
$t=1$.  Then we can  use Lemma 14 in \cite{kn:D3}, much as in Lemma 8,  to get
\begin{equation} \sqrt{J(1)} \leq c \ \max( (u'(2)-u'(\frac{1}{2})),  \Emax). \end{equation}
 
Now
$ J(1)\geq B^{-2}  J(2) \exp( u'(2)-u'(1))= B^{-2} J(2) \exp(\delta_{0})$ by Proposition 12 and so, using our lower bound
of Lemma 10,
$J(1) \geq c \exp(\delta_{0})$. If $\exp(\delta_{0})$ is large compared with
$\Emax^{2}$ we must have
$ \sqrt{J(1)} \leq  c (u'(2)- u'(\frac{1}{2})= c(\delta_{0}+ \delta_{1}),
$
so we get
\begin{equation} e^{\delta_{0}} \leq c (\delta_{0}+\delta_{1}). \end{equation}
Thus $\delta_{1} \geq f(\delta_{0})$ where $f$ is the function
$$  f(\delta)= c^{-1} e^{\delta} - \delta. $$
We can obviously choose a  $\udelta>1$ such that if $\delta\geq \udelta$ we have
$$f(\delta)\geq \delta^{2}\geq \delta\geq \udelta.$$
Then if $\delta_{0}\geq \udelta$ we have $\delta_{1}\geq \delta_{0}^{2}\geq
\delta_{0}$.
Now the whole set-up is invariant under rescaling by a power of $2$, so we
also have  $\delta_{n+1} \geq \delta_{n}^{2}$. hence $\delta_{n}\geq \delta_{0}^{2^{n}}$. But by an easy argument this would imply that $u(t)$ is unbounded as $t\rightarrow
0$, contrary to what we know. So we deduce that in fact either $\delta_{0}\leq
\udelta$ or $\exp(\delta_{0})\leq c \Emax^{2}$.  Now the
statement for all $n$ follows by rescaling.

\

Now set 
\begin{equation} \Delta(t)= t^{2} \max \{ J(x_{1}, x_{2}) : x_{1}+ x_{2}= 2t; \vert x_{1}-
x_{2}\vert \leq t/10\}.\end{equation}
Note that the factor $t^{2}$ in the definition makes this invariant under rescaling.

We introduce a parameter $\epsilon \in (0,1)$, to be fixed later, and consider
the function
   \begin{equation}  F_{\epsilon}(t)= E(t) + \epsilon \Delta(t). \end{equation}
   After scaling we can suppose this achieves its maximal value $\Fmax$ at
   $t=1$, we write $E=E(1), \Delta=\Delta(1)$.
   Now using the bound from Proposition 13, and (44) we get
   $\sqrt{\Delta}\leq c\max(E, \log \Emax)$ so
   $$\Emax\leq \Fmax\leq E + (E+ \log \Emax)^{2} . $$
   This gives
   $$ \Emax\leq c( E^{2} + (\log \Emax)^{2}). $$
   Thus \begin{equation}\Emax\leq c E^{2}. \end{equation} We can suppose that $E$ is large (for otherwise
   $\Emax$ is not too large) then we get
   \begin{equation} \delta_{n} \leq 4 \log E. \end{equation}
   
   \
   
   \

    We normalise $u$, under the addition of affine-linear functions, at the
    point $(1,1)$. Then summing the $\delta_{n}$, using the bound (49) and
    integrating the resulting bound on the $\frac{\partial u}{\partial t}$
    we see that
    the variation of $u$ over compact subsets of the diagonal $\{s=0\}$ is $O(\log E)$, which is small compared with the variation
across the orthogonal line $\{t=1\}$, since the latter is at least $E$, by definition.
    More generally we have
    \begin{lem}
    For any $t_{2}>1$ and $\sigma$ with $\vert\sigma \vert\leq 1/2 $ the variation of $u$ on the intersection of the
    line $\{ x_{1}- x_{2}= 2\sigma \}$ with the triangle $\{x_{1}+x_{2}\leq 2t_{2}\}$ is bounded by $c\log E$, where $c$ depending only on  $t_{2}$.
    \end{lem}
    
    We know that $u$ is $O(\log E)$ on the diagonal and it follows from the
    definitions that 
    $ u $ is $O(\Emax)$ on the triangle $\{x_{1}+x_{2}\leq 3t_{2}\}$,
    say. This means that the size of the derivative of $u$ is $O(\Emax d^{-1})$
    where $d$ is the distance to the boundary. Then by applying Theorem 5 of \cite{kn:D3}  we deduce that
   \begin{equation}J\leq c \Emax^{2} d^{-4}.\end{equation}   Consider the line $\{ x_{1}-x_{2}= 2\sigma\}$, where we can suppose $\sigma\geq
  0$, and parametrise this line by $x_{1}= 2\sigma+\tau, x_{2}=\tau$. By
  applying Proposition 12 and the lower bound of Lemma 10 we see that 
    $$    \vert \frac{\partial u}{\partial \tau} \vert \leq c \log( c E^{4} \tau^{-4}),
    $$ where we have used (48) to replace $\Emax$ by $E$.  
    Integrating this we obtain the result.

    \

     In the next two subsections we prove the following two propositions.
    
  \begin{prop}
  There is a $k_{0}$, independent of $\epsilon$, and a function $\mu(\epsilon)$
  such that if at a interior maximum point for $F$ we have $\Delta\geq k_{0}
  M$ then $\Emax\leq \mu(\epsilon)$.
  \end{prop}
  \begin{prop}
  For any $k$ there is an $\epsilon(k)$ and $\nu(k,\epsilon)$ such that if $\epsilon\leq \epsilon(k)$
  and if at an interior
  maximum point for $F$ we have $\Delta\leq k M$ then $\Emax\leq \nu(k,\epsilon)$.
  \end{prop}
  These two propositions complete the proof of Proposition 11.
   For we fix $\epsilon= \epsilon(k_{0})$ and then
   at an interior maximum we have 
    $$\Emax\leq \max( \mu(\epsilon(k)), \nu(k_{0}, \epsilon(k_{0})) ). $$

    We will use a simple principle in the proofs of both of these Propositions.
    Write $\xi_{s}, \xi_{t}$ for the partial derivatives of $u$ with respect
    to the variables $s,t$.
    Given a point $p$ and real numbers $\alpha, \beta_{1}, \beta_{2}$  with
    $\beta_{1}< \beta_{2}$,  let $S=S(p;\alpha, \beta_{1}, \beta_{2} )$ be the set of points $(x_{1}, x_{2})$ where
  \begin{equation}    \beta_{1}\leq \xi_{s}(\ux) \leq \beta_{2}\ ,\ \xi_{t}\leq \alpha+ \xi_{t}(p).
  \end{equation}
  
  \begin{lem} We have
     $$  B^{-2} (\beta_{2}-\beta_{1}) J(p)^{-1} \leq \Area(S)\leq B^{2} (\beta_{2}-\beta_{1}) J(p)^{-1}.
     $$
     \end{lem}  
     For the area of $S$ is 
  \begin{equation} \Area(S)= \int_{\Pi} J^{-1} d\xi_{s} d\xi_{t},  \end{equation}
  where $\Pi$ is the region in the $(\xi_{s}, \xi_{t})$ plane defined by
  the inequalities (51), and we have abused notation by regarding $J$ as a function
  of $\xi_{s}, \xi_{t}$ in the obvious way. Now the volume bound of Proposition
  12 gives
  \begin{equation} B^{-2} J(p) e^{\xi_{t}(p)-\xi_{t}} \leq J(\xi_{s}, \xi_{t}) \leq B^{2} J(p) e^{\xi_{t}(p)-\xi_{t}}, \end{equation}
and the result follows by integrating the exponential function over $\Pi$.

  \subsection{Proof of Proposition 14}

    We fix values $t_{0}, t_{1}, t_{2}$, say for definiteness $t_{0}=1/5,
    t_{1}=1/4$ and $t_{2}=2$.
    Let $R$
   be the rectangle $\{ \vert s\vert \leq 1/10, t_{0} \leq t \leq t_{2} \}$.
   
    Recall that the definition
   of $\Delta$ involves maximising over an interval $\vert s\vert \leq t/20$.
   Suppose that the maximum is achieved at a point $p$, where $t=1$ and $s=s_{0}$.
   (Of course, we can suppose $t=1$ by rescaling.)
   So $\vert s_{0}\vert\leq 1/20$ and $p$ lies inside $R$.
   
   Now the proof proceeds by the following steps.

   \
   
   {\bf Step 1} {\it Claim:  There is a $c_{1}$ such that $\vert \frac{\partial u}{\partial s} \vert
   \leq c_{1} E$ on $R$.}

\

For on the line segment $\{ t=1, \vert s\vert \leq 1\}$ we have a bound $\vert
u\vert \leq
 E \vert s\vert$. Using Lemma 11, this gives an $O(E)$ bound on $u$ over the interior
 region $\vert s\vert \leq 1, t\leq 2$.  Since $R$ lies within the interior
 of this set, convexity gives an $O(E)$ bound on the derivative over $R$.

 \

 Now we consider the set $S=S(p; \alpha, -c_{1} E, c_{1} E )$, with $c_{1}$ as above. The curve $\{ \frac{\partial
 u}{\partial t}=\alpha\}$ is the graph of a function $t=\tau(s)$. By item
 (1) above the intersection $S\cap R$ is just the set defined by the three
 conditions
$$-1/10\leq s\leq 1/10, t_{0}\leq t\leq t_{2}, t\leq \tau(s).$$

\

\

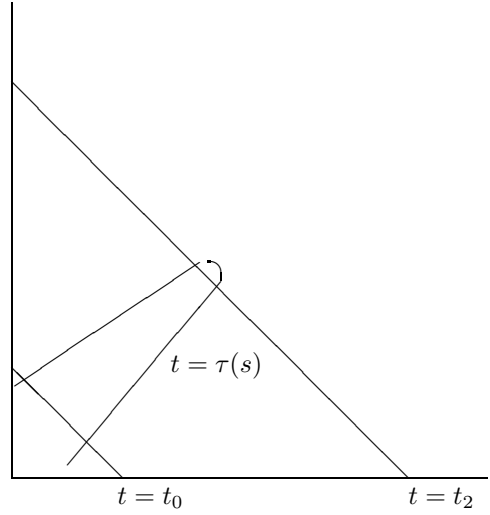
 \begin{picture}(180,180)(-30,-30)
\put(0,0){\line(1,0){180}}
\put(0,0){\line(0,1){180}}
\put(150,0){\line(-1,1){150}}
\put(42,0){\line(-1,1){42}}
\put(150,-10){$t=t_{2}$}
\put(40,-10){$t=t_{0}$}
\put(60,40){$t=\tau(s)$}
\put(74,75){\oval(10,14)[tr]}
\put(21,5){\line(5,6){58}}
\put(1,35){\line(3,2){70}}
\end{picture} 

\

\

{\bf Step 2}
 \ \ \ {\it Claim:  We can choose $\alpha>0$, depending only on $B$, so that for any point
 $q$ on the graph $t=\tau(s)$ and any point $p'$ with $s=s_{0}, t\leq 1$ we have $$J(q)< \frac{t_{2}- t_{1}}{1-t_{1}} J(p'). $$}

For, since the partial derivative $\frac{\partial u}{\partial t} $ is monotone
on the line $s=s_{0}$ we have
$$     J(p')\geq B^{-1} J(p), $$
whereas, by the inequality of Proposition 12, 
$$   J(q)\leq B e^{-\alpha} J(p). $$
So we just need to choose $\alpha>0 $ and bigger than $\log\left( B^{2}
 \frac{t_{2}-t_{1}}{1-t_{1}}\right)$.

Now we fix $\alpha$ as above. By Lemma 12, the area of $S$ is at most $ c E/J(p)=
c E/\Delta$. So, by choosing $k_{0}$ large (as allowed in the statement of Proposition
14), we can suppose the area of $S$
is as small as we please. Fix a  suitably small number $\delta$--- for
 definiteness we can take $\delta=1/100$---and choose $k_{0}$ so that the area of $S$
is less than $\delta (t_{1}- t_{0})$.

Note that, since $\alpha>0$, we have $\tau(s_{0})>1$, by monotonicity of the partial
derivative.

\

\
 
{\bf Step 3}\ \ \ {\it Claim: There are $s_{-}, s_{+}$ with $\vert s_{\pm}- s_{0}\vert \leq \delta$
and $s_{-}< s_{0}< s_{+}$ such that $\tau(s_{\pm})\leq t_{1}$.}

\

If there is no such $s_{+}$ then $S$ contains the rectangle $t_{0}\leq t
\leq t_{1}, s_{0}\leq s\leq s_{0}+\delta$. (Notice that our choices imply
that this rectangle lies inside $R$.) But this contradicts the fact that
the area of $S$ is less than $\delta(t_{1}- t_{0})$. Similarly for $s_{-}$.

\

To sum up so far we have shown that the set $S$ must contain a very thin
\lq\lq finger'', extending out from the region $\{ t\leq t_{1}\}$ and containing
the point $p$ where $t=1$.

\

Let $s_{+}$ be the least among the values satisfying the conditions of the
claim above
and $s_{-}$ be the largest. Then $\tau(s_{\pm})=t_{1}$ and $\tau>t_{1}$
on the open interval $(s_{-}, s_{+})$.

Let $\Omega$ be the set where $s_{-}\leq s\leq s_{+}$, $t_{1}\leq t\leq t_{2}$
and $t\leq \tau(s)$. For $t'\in [t_{1}, t_{2}]$ let $I_{t'}$ be the intersection of
$\Omega$ with the line $t=t'$ and let $j(t')$ be the maximum of $J$ over
$I_{t'}$. Thus $j(1)\geq \Delta$.  Let $G\subset [t_{1}, t_{2}]$ be the set of values $t^{*}$ such that for all $t'\in (t_{1}, t^{*}]$ the maximum $j(t')$
is attained at an interior point of  $I_{t'}$ (i.e. not at points in
the graph of $\tau$).

\

\

{\bf Step 4} \ \ {\it Claim: $1$  is contained in $G$.}

\

This follows from the Claim in Step 2, since for $t'\leq 1$ the point $p'$ with co-ordinates
$t=t', s=s_{0}$ lies in
$I_{t'}$ and $J(p')$ is strictly less than the value of $J$ at any point
on the graph.

\

Now the crucial idea in the proof is to show that this thin \lq\lq finger''
must actually extend to meet the line $\{t=t_{2}\}$.

\

\

{\bf Step 5}\ \  {\it Claim:  If  $t^{*}\geq 1 $ and $t^{*}\in G$
  then $j^{-1}$ is a concave function on the interval $[t_{1}, t^{*}]$.
}

\
This is similar to the proof of Proposition 9.
The function $F=J^{-1}$ satisfies the linear equation $u_{ij} F^{ij}=
-A$ and $A\geq 0$. Then the assertion follows from the maximum principle
applied to $F-c t$ for suitable values of $c$. 

\

{\bf Step 6}\ \  {\it Claim:  $t_{2}$ is in  $G$. }

\

This follows from a continuity argument. From its definition, $G$ is open. So long as $t^{*}$ lies
in $G$ we have  $$j(t^{*})^{-1} \leq \frac{t^{*}-t_{1}}{1-t_{1}} j(1)^{-1}
\leq  \frac{t^{*}-t_{1}}{1-t_{1}} J(p)^{-1}, $$
by convexity. Suppose $1\leq t^{*}\leq  t_{2}$. Recall that we arranged that for any point $q$ on the graph of $\tau$
$$J(q)^{-1}> \frac{t_{2}-t_{1}}{1-t_{1}} J(p)^{-1}.$$
The strict inequality implies that $G$ is closed.

\

\

{\bf Step 7}\ \  {\it Claim: There is a point $p''$ with co-ordinates $s''\in (s_{-}, s_{+})$  and $t=t_{2}$ such that $J(p'')>
 \frac{1-t_{1}}{t_{2}-t_{1}}J(p)$}
 
 This follows from the concavity of $j^{-1}$, as above.
 
 Now by the choice of $\delta$ we have $\vert s''\vert \leq t_{2}/20$ so the
 point above is one of those considered in the definition of $\Delta(t_{2})$
 and we have
 $$  \Delta(t_{2}) \geq \frac{1-t_{1}}{t_{2}-t_{1}} t_{2}^{2} \Delta. $$
 
 Now with the definite choices of $t_{i}$ made above this inequality  is
 $\Delta(t_{2}) \geq (1+\sigma) \Delta$ with $\sigma=5/7>0$.
 
 \
 
 \

  We can now complete the proof. From Lemma 11 we know that $u(0,0)$ is $O(\log E)$ and the convexity of $u$ on the boundary implies
 that $$u(2t_{2}, 0)\geq t_{2} u(1,0) - c \log E, u(0,2t_{2})\geq u(0,1) -c
 \log E.$$
 Also $u(t_{2}, t_{2})$ is $O(\log E)$, again by Lemma 11, so from the definition
 of $E(t)$ we have
 $$ E(t_{2}) \geq E - c \log E. $$
 So $$F_{\epsilon}(t_{2}) \geq E + \epsilon (1+\sigma) \Delta- c \log E\geq
 F_{\epsilon}(1) +\epsilon \sigma
 k_{0} E - c\log E . $$
 Proposition 14 follows from (48) and the fact that $F_{\epsilon}(t_{2})\leq F_{\epsilon}(1)$.

 \subsection{Proof of Proposition 15}

 \

 Recall from the statement of the Proposition that we are supposing that $\Delta \leq k E$. The main idea in the proof will be in part complementary
to that of Proposition 14, in that we invoke a {\it lower} bound on the area
of a suitable
set $S$.  As before we suppose that the maximum of $F_{\epsilon}$ is attained
at $t=1$, and let $p$ be the point on the line $\{ t=1\}$ where the maximum
in the definition of $\Delta$ is achieved.  Our argument again employs certain parameter values $t_{0}, t_{2}$
for the $t$-coordinate, but this time we will choose $t_{0}<1$ very small, so that
\begin{equation}   \frac{1}{2} t_{0}^{2} \leq \frac{1}{100 B^{2}k}, \end{equation}
and $t_{2}>1$ very large, so that
\begin{equation}    t_{2}^{2} \geq 100 B^{4} k . \end{equation}

For any $t$ we have
$$  E(t)\leq F_{\epsilon}(t) \leq F_{\epsilon}(1)=E+\epsilon \Delta \leq
(1+\epsilon k) E. $$ 
Write $U(t)=u(2t,0)+u(0,2t)$, so $E=U(1)$.  Over the fixed range,  $t\leq t_{2}$, our
bound in Lemma 11, on the diagonal, gives $U(t)\leq t E(t)+ c \log E$. Hence
\begin{equation} U(t) \leq t (1+\epsilon k) E + c\log E = t(1+\epsilon k)
U(1) + c\log E\end{equation}

 When $\epsilon$ is small, convexity of the function $U$ forces $E^{-1} U(t)$ to be close to the linear function $t$ (assuming of course that $E$ is large), over the
range $t\leq 1$. Further, each summand $u(2t,0), u(0,2t)$ is positive and
convex and this forces
$$   u(2t,0) = t u(2,0) + O(\log E+  \epsilon E)\ \ , \ \ u(0,2t)= t u(0,2) + O(\log E+\epsilon E). $$
To express this differently, write $u(2,0)= \lambda_{1}E, u(0,2)=\lambda_{2}E$,
so $\lambda_{i} \geq 0$ and $\lambda_{1}+\lambda_{2}=1$. Define
$$ V(x_{1},x_{2})= \frac{1}{2} \max\left( \lambda_{1}(x_{1}-x_{2}), \lambda_{2}(x_{2}-x_{1}\right).
$$ 
Then on compact subsets of the boundary of the quarter-plane the function
$E^{-1} u$ differs from $V$ by $O( E^{-1} \log E+\epsilon)$. But now Lemma 11 implies
that $E^{-1} u$ differs from $V$ by $O(E^{-1} \log E+\epsilon)$ over the whole region
$ \{t \leq t_{2}\}$ in the quarter plane.
Set $$\beta_{1}= \left( \frac{1}{2} (\lambda_{1}-\lambda_{2})-\frac{1}{4}\right)
E\ \ , \ \ \beta_{2} =\left(\frac{1}{2}(\lambda_{1}-\lambda_{2})+\frac{1}{4}
\right) E. $$
(The reader may find it easiest to think first of the symmetrical case
when $\lambda_{1}=\lambda_{2}=1/2 $.) We consider the set $S=S(p;\alpha,\beta_{1},
\beta_{2})$, where $\alpha>0$ and
\begin{equation}    e^{\alpha} \geq 10 B^{2} k. \end{equation}
(The reason for these choices will emerge presently.) The crucial observation is
\begin{lem}
Thee intersection of $S$ with the set $\{ t_{0}\leq t\leq t_{2}\}$ is
contained in a strip $\vert s \vert \leq \eta$ where $\eta= c  \epsilon$,
once $E$ is sufficiently large.
\end{lem}
The proof is straightforward, using the preceding discussion. (Only the constraint
$\beta_{1}\leq \xi_{s} \leq \beta_{2}$ is relevant here: the statement
is valid for any $\alpha$.)

\

Now we proceed with the following steps.

\

{\bf Step 1}\ \  
{\it  Claim: If $E$ is sufficiently large, the variation of $\frac{\partial u}{\partial s}$ on the line $t=1$
across the strip $\vert s \vert \leq \eta$ is at least $E/2$. }

We know that $\frac{\partial u}{\partial s}$ varies from $-\infty$ to $\infty$
across the whole interval $t=1, \vert s \vert \leq 1$ and by the Lemma 13 the points
where values $\beta_{1}, \beta_{2}$ are attained must lie in this strip.
Then the claim follows from the fact that $\beta_{2}-\beta_{1}= E/2$.

\

\

{\bf Step 2}\ \  {\it Claim: There is a $k'\geq 1/10$,
such that, if $E$ is sufficiently large,  there is a point $p'$ in the segment
$\vert s \vert \leq \eta, t=1$  with $\beta_{1}\leq \xi_{s}(p')\leq \beta_{2}$  and  $J(p') \geq k'E$.}
 
To see this we use an integral identity just as in \cite{kn:D3}, Lemma 17. This gives a formula, in terms of $A$, for the integral over the line segment
$t=1, \vert s\vert \leq 1$ of $u^{tt}$ (in an obvious notation). Using the
formula for the inverse of a $2\times 2$ matrix, we can write this is
$$  \int J^{-1}  d \xi_{s},  $$
where $\xi_{s}=\frac{\partial u}{\partial s}$ is regarded as a parameter on the line segment $\vert s\vert \leq
1, t=1$. From this formula one sees that the integral is at most $5$, when
$\vert A \vert_{L^{\infty}}\leq 1$, as we are supposing. In particular the same integral over the sub-segment $\vert s\vert
\leq \eta, t=1$ is bounded above by $k'/2\leq 5$ say. 
 Since the variation in $\xi_{s}$ over this subsegment at least $E/2$ there
 must be a point where  $J\geq k'E$, with $k'\geq 1/10$.

\

\

{\bf Step 3}\ \  {\it Claim: If $\epsilon$ is sufficiently small and $E$ is sufficiently
large then  $\Delta\geq k'E$. } 
We just  choose $\epsilon$ so that $\eta<1/20$ and the point $p'$ is one
of those considered in the definition of $\Delta$.

\

{\bf Step 4} {\it Claim: If $\epsilon$ is sufficiently small, and $E$ sufficiently
large, the set $S$ intersects the line segment $\{ t=t_{2}, \vert
s \vert \leq  \eta\}$} 

\

First, since $\frac{J(p')}{J(p)} \geq (10 k)^{-1}$ it follows from the choice
of $\alpha$ and Proposition 12 that $p'$ lies in  $S$. By Lemma 12, the area of $S$ is at least
$\frac{E}{2\Delta B^{2}}\geq (2 k B^{2})^{-1}$ and by the choice of $t_{0}$ this is more than
twice the area of the triangle $\{t\leq t_{0}\}$. Suppose $S$ does not intersect
the line segment as claimed. Since $p'$ lies in $S$ and $S$ is connected
it follows from Lemma 13 that $S$ is contained in the union of the triangle $\{
t\leq t_{0}\}$ and the strip $\{ t_{0}\leq t \leq t_{2}, \vert s\vert \leq
\eta \}$. So the area of this strip must be at least half the area of $S$.
But the area of strip is $2(t_{2}- t_{0})\eta= 2 c (t_{2}-t_{0}) \epsilon$,
so this is impossible when $\epsilon$ is sufficiently small.

\

\

{\bf Step  5} \ {\it Claim: $\Delta(t_{2})\geq 2 \Delta$}.

Consider the point $p''$ whose existence is established in the previous step.
Using Proposition 12, $J(p'')\geq B^{2} e^{-\alpha}J(p)= B^{2} e^{_\alpha} \Delta$. We
can assume that $\eta$ is small, so $p''$ is a point considered in the definition
of $\Delta(t_{2})$ and $\Delta(t_{2}) \geq t_{2}^{2} J(p'')$. Now the claim
follows from the choice of $t_{2}$.

\

Now we can complete the proof. By convexity of the function $U$ and fact
that $U(0)$ is $O(\log E)$ we have
$$  U(t_{2}) \geq t_{2} E - c \log E. $$
From the bound on $u(t_{2}, t_{2})$ we deduce that
$$   E(t_{2}) \geq E - c  \log E. $$
So $$F_{\epsilon}(t_{2}) \geq E+ \epsilon \Delta(t_{2})- c\log E\geq E+ 2
\epsilon \Delta- c\log E. $$
Now the fact that $F_{\epsilon}(t_{2})\leq F_{\epsilon}(1)$ gives
$$   E+ 2\epsilon \Delta- c\log E\leq E + \epsilon \Delta, $$
so
$   \epsilon \Delta \leq c\log E $. Now by Step 3, $\Delta\geq k' E$ so
$$ \epsilon k' E \leq c \log E, $$
which gives the required bound on $E$.

\subsection{Completion of proof of Main Theorem}

We need to control a solution $u^{(\alpha)}$ in a neighbourhood of a vertex, and we can
take standard co-ordinates around the vertex as in the previous section.
For each point $p=(p_{1}, p_{2})$ we define $D_{1}(p)$ as in (3.2)
$$   D_{1}(p)= \frac{1}{ p_{1}} ( u(0,p_{2}) - u(p_{1}, p_{2})- \frac{\partial
u}{\partial x_{1}}(p_{1}, p_{2}). $$
Of course we have a similar quantity $D_{2}(p)$ defined by interchanging the co-ordinates. We prove
\begin{prop}
There is a {\it a priori} bound $D_{1}(p), D_{2}(p) \leq D$, valid for all
solutions $u=u^{(\alpha)}$ in our sequence, and all points $p$ near to a
vertex.
\end{prop}
Given this it is straightforward to adapt the proofs of proposition 6 to deduce Theorem
1, arguing just as in (4.1).

When $p_{1}=p_{2}$ a bound on $D_{i}(p)$ follows immediately from what we
have proved in the previous section.  More generally, the only new issues
arise when $p_{1}$ is much less than $p_{2}$. We adapt the argument of Section
3. We choose a point where $D_{1}(p)$ is maximal and by rescaling we can suppose
that $p_{2}=1$. If $u$ is normalised at $(1,1)$, as in the previous subsection,
then we have {\it a priori} $L^{\infty}$ bounds on $u$ over compact subsets.
The only difficulty in applying the argument of Section 3 would occur if the \lq\lq
scale''
$\lambda$ is not small. Then the limit $\ufl$ would still be defined only
on a quarter plane and the problem would come when, taking the limit in seeking
a contradiction, the $\ufl$ converge to an affine-linear function on the boundary. To rule this out we need an {\it a priori} \lq\lq strict convexity''
bound on the restriction of $u$ to the $x_{2}$-axis. Thus the crucial thing
is to prove
\begin{prop}
There are $r>1$ and $\eta>0$ such that 
$$  \frac{\partial u}{\partial x_{2}}(0,r)- \frac{\partial u}{\partial x_{2}}(0,1)\geq
\eta, $$
for all functions $u$ obtained by rescaling a $u^{(\alpha)}$
\end{prop}
Given this proposition it is very to adapt the arguments of Section 3 to
prove Proposition 16, on the lines indicated above.

\

To prove Proposition 18, we begin by considering the derivative $\xi_{t}=\frac{\partial
u}{\partial t}$ on the diagonal. We have upper and lower bounds on the determinant
$J$ at the point $(1,1)$ and, under rescaling,  these give
$$  c' t^{-2}\geq   J(t,t)\geq c t^{-2}. $$
Then Proposition 12 gives $$   \xi_{t} \geq 2 \log t - c $$
for $t\geq 1$. Now write $\frac{\partial u}{\partial x_{2}}= \xi_{t}-\xi_{s}$
so
$$    \frac{\partial u}{\partial x_{2}} (t,t) \geq 2 \log t +\xi_{s}(t,t)
- c.$$

The $L^{\infty}$ bound on $u$ implies a bound on $\frac{\partial u}{\partial
x_{2}}(1,1)- \frac{\partial u}{\partial x_{2}}(0,1)$. By considering the
rescaling behaviour we get a fixed bound on $\frac{\partial u}{\partial x_{2}}(t,t)-\frac{\partial
u}{\partial x_{2}}(0,t)$ for all $t$. So
\begin{equation} \frac{\partial u}{\partial x_{2}}(0,t)\geq 2\log t + \vert  \xi_{s}(t,t)\vert - c. \end{equation}
Thus it suffices to show that $\vert  \xi_{s}(t,t)\vert$ is small, for large
$t$, compared with $2 \log t$. To this end we first define
$$   T= T(u)=\vert \xi_{s}(2,2)-\xi_{s}(1,1) \vert. $$
 (In fact $u$ has been normalised so that
$\xi_{s}(1,1)=0$, but it clearer to write the definition this way.) By our
$L^{\infty}$ bounds we have $T\leq C_{0} $ say.

For integer $\mu>1$ let $\Omega_{\mu}$ be the region $\{ 2^{-\mu} \leq t \leq 4\}$
in the quarter-plane. Let $E_{\mu}(u)$ be the integral of the quantity $\vert
F\vert^{2}$ (as defined in \cite{kn:D2}, \cite{kn:D3}) over $\Omega_{\mu}$. (This is essentially the square of the $L^{2}$ norm of the Riemann curvature tensor over the corresponding piece
of a $4$-manifold.) Given a positive number $C$, let ${\cal A}_{C}$ be the
set of convex functions $u$ on the closed triangle $\{t\leq 4\}$, normalised
at $(1,1)$ such that
\begin{enumerate}
\item $u$ satisfies equation (1), with $\Vert A\Vert_{C^{2}}\leq C$;
\item $\Vert u\Vert_{L^{\infty}}\leq C$;
\item $\det(u_{ij})\geq C^{-1}$ everywhere;
\item $\Vert V \Vert_{L^{\infty}}\leq C$, where $V$ is the vector field 
associated to $u$;
\item $u$ satisfies Guillemin boundary conditions, with the standard measure,
along the $x_{i}$-axes;
\item $\det(u_{ij}) \geq C^{-1} t^{-2}$ at the point $(\frac{t}{2}, \frac{t}{2})$.
\end{enumerate}
Then we have
\begin{lem}
For any $C, \epsilon>0$ there is an integer $\mu$ and a $\delta>0$ such that if $u\in {\cal A}_{C}$ and
$E_{\mu}\leq \delta$ we have $T(u)\leq \epsilon$.
\end{lem}
For our application we fix $\epsilon<2\log 2$. 

Assuming this lemma for the moment, we complete the proof of Proposition 17. For integers $n\geq 1$ set
$$  T_{n} = \vert \xi_{s}(2^{n}, 2^{n})- \xi_{s}(2^{n-1}, 2^{n-1}) \vert.
$$
Thus, rescaling by a factor $2^{n}$, we can apply our bounds on $T$ to give
bounds on $T_{n}$. For suitable $C$, all of the conditions defining ${\cal A}_{C}$ hold for our $u^{(\alpha)}$,  so $T_{n}\leq C_{0}$ for all $n$, 
and either $T_{n} \leq
\epsilon$ or the integral of $\vert F\vert^{2}$ over the region
$\{     2^{n-\lambda}\leq t \leq 2^{n+1}\}$ exceeds $\delta$. Now use the fact
that we have a fixed bound on the $L^{2}$ norm of $F$ over the whole polygon,
so we can find a large integer $M$ such that 

$$  \int_{P}\vert F\vert^{2} \leq M \delta. $$
It follows that there are at most $M (\lambda+1)$ values of $n$ for which $T_{n}$ exceeds
$\epsilon$. Thus
\begin{equation}  \vert \xi_{s}(2^{n}, 2^{n})-\xi_{2}(1,1)\vert \leq M (\lambda+1) C_{0}+ n \epsilon.\end{equation}
Combining (59) (with $t=2^{n}$ ) with (58) we establish Proposition 17.

It only remains to prove Lemma 14. Arguing by contradiction, we suppose
 we have a sequence of  functions $u^{(\beta)}\in {\cal A}_{C}$ with the integral of $\vert
F\vert^{2}$ over $\Omega_{p}$ tending to zero, for a sequence $p$ tending
to infinity, and the sequence of quantities $T(u^{(\beta)})$ does not tend to zero. Using the first three conditions in the definition of ${\cal A}_{C}$ and the arguments of \cite{kn:D2} we see
that, taking a subsequence, we can suppose the sequence converges in $C^{4}$
on compact subsets of the interior. The limit clearly has $F=0$, i.e. describes
a flat metric. The vector field $V$ associated to
this limit is constant and  a straightforward Stokes' Theorem argument (similar
to that in Proposition 5 ), using the boundary conditions and the fourth condition before taking the limit, shows that $V_{\infty}$ is the vector
field $\frac{\partial}{\partial x_{1}}+ \frac{\partial}{\partial x_{2}}$.
There is a simple classification of locally-defined functions $u_{\infty}$
with $F=0$ and it is easy to read off from this that in the limit $\xi_{s}$
is constant along the diagonal, which gives the desired contradiction. (The
point here is that the fifth condition forces the limit to blow up at the
origin.)

\section{Blow-up limits}

\subsection{The Joyce construction}

We recall a construction, due to Joyce, of explicit solutions of equation
(1), with $A=0$, that is, metrics of zero scalar curvature. The original reference
is \cite{kn:J}, but we follow the approach of Calderbank and Pedersen in \cite{kn:CP}. An elementary
derivation of this construction (and a generalisation to other equations)
from the point of view of this paper is given in the note \cite{kn:D4}. 

\

Consider the linear PDE for a function $\xi(r,H)$, where $r>0$, 

\begin{equation} \frac{\partial^{2} \xi}{\partial H^{2}} + r^{-1} \frac{\partial}{\partial
r}\left( r \frac{\partial \xi}{\partial r}\right)= 0. \end{equation}
This is familiar as the equation defining axi-symmetric harmonic functions
in cylindrical co-ordinates on $\bR^{3}$.
 Given a pair of solutions $\xi_{1}, \xi_{2}$ to (60), we set
  $$  P_{i}= \frac{\partial \xi_{i}}{\partial H} \ , \ Q_{i}= \frac{\partial
  \xi_{i}}{\partial r}. $$
  and write $\Delta= P_{1}Q_{2}-Q_{1}P_{2}$. We assume that $\Delta>0$ everywhere. We introduce two further angular co-ordinates $\theta_{1}, \theta_{2}$ and
consider the four dimensional Riemannian metric
  $$  g= \frac{r\Delta}{2} (dH^{2} + dr^{2}) + \frac{r}{2\Delta}\left( (P_{2}^{2}
  + Q_{2}^{2}) d\theta_{1}^{2} - 2(Q_{1} Q_{2}+ P_{1}P_{2})d\theta_{1}d\theta_{2}+
  (P_{1}^{2}+ Q_{1}^{2}) d\theta_{2}^{2}\right). $$
  The main result is that this is a Kahler metric of zero scalar curvature.
  To relate this to the equation (1), we introduce another linear equation
   \begin{equation} \frac{\partial^{2} x}{\partial H^{2}} + r \frac{\partial}{\partial
r}\left( r^{-1} \frac{\partial x}{\partial r}\right)= 0. \end{equation}
 Given a solution $\xi(r,h)$ to (60), the first order system
 $$   \frac{\partial x}{\partial r}= r \frac{\partial \xi}{\partial H}\ \
 ,\ \  \frac{\partial
 x}{\partial H}= - r \frac{\partial \xi}{\partial r}, $$
 is consistent and has a solution $x$, unique up to a constant. Further $x$
 satisfies the equation (61).  So
 starting with a pair of solutions $\xi_{1}, \xi_{2}$ to (60) we get a pair
 of solutions $x_{1}, x_{2}$ to (61), but we introduce a sign and interchange
 of labels so that
 $$  \frac{\partial x_{1}}{\partial r}= \frac{\partial \xi_{2}}{\partial
 H}\ \ , \ \ \frac{\partial x_{2}}{\partial r} = - \frac{\partial \xi_{1}}{\partial
 H}, $$
 {\it etc.}. Then $x_{i}$ and $\xi_{j}$ give the co-ordinates we are familiar
 with in this paper: the equations imply that the $1$-form $\sum \xi_{i} dx_{i}$
 is closed, so there is, at least locally, a function $u$ with
 $ du= \sum \xi_{i} dx_{i}$. If we regard $u$ as a function of $(x_{1}, x_{2})$
 then we get a solution of (1), with $A=0$. Conversely any solution arises
 in this manner away from  the critical points of $\det(u_{ij})$.
 In fact the construction gives $r=\det(u_{ij})^{-1/2}$.

  An important special case occurs when $\log J$ is an affine-linear function of $\xi_{1}, \xi_{2}$.  In differential
geometric terms, our metric is then Ricci-flat. Making an affine change
of variable we may suppose that $ \xi_{2}=\log r$.
So $P_{2}=0$ and $Q_{2}=r^{-1}$. It is easy to check that the metric $g$
is the same as that
given by the well-known  Gibbons-Hawking construction, using the harmonic function
$\frac{\partial \xi_{1}}{\partial H}$ on $\bR^{3}$.

 With  this background in place, we can move on to consider the particular
 metrics we are interested in. Consider first the case of flat space, so
 $$u= x_{1} \log x_{1} + x_{2} \log x_{2} $$ 
 on the quadrant $\{x_{1}, x_{2}>0\}$. Then one finds that
 $$\xi_{1}= \log F_{-}(H,r), \xi_{2}= \log F_{+}(H,r)$$ 
 where
 $$  F_{\pm}(H,r) =\frac{1}{2}\left( \pm H +\sqrt{ H^{2} + r^{2}}\right). $$
 The harmonic function $\log F_{-}$ is the potential associated to a uniform charge
 distribution on the half-line $r=0, H>0$ and $\log F_{+}$ to the half-line $r=0,H<0$.
 So, as we see from the formulae,  $F_{\pm}$ has a logarithmic singularity
 along the corresponding half-line and
 $$  F_{+}+ F_{-}= 2 \log r. $$
 The corresponding functions $x_{i}$ are just $x_{1}= F_{-}, x_{2}= F_{+}$.

 Now  given $a_{1}, a_{2}>0$ set
 $$  \xi_{1}= \log F_{-}(r,H)- a_{2} H +1, \xi_{2}=\log F_{+}(H,r)+a_{1} H+1. $$
(The addition of the constant $1$ makes no change to the geometry but will
be convenient later.)
These are obviously harmonic functions and the corresponding functions $x_{i}$
are
$$  x_{1}= F_{-}+ \frac{a_{1} r^{2}}{4}\ , \ x_{2}= F_{+} + \frac{a_{2} r^{2}}{4}.
$$
We want to find the \lq\lq symplectic potential'' $u(x_{1}, x_{2})$ which describes
this solution. Set $y_{1}= F_{-}, y_{2}= F_{+}$ so that
$y_{2}-y_{1}= H$ and $y_{1}y_{2}= r^{2}/4$. Thus
$$  x_{1}= y_{1}+ a_{1} y_{1} y_{2}\ , x_{2}= y_{2}+ a_{2} y_{1} y_{2}. $$
and $$\xi_{1}= \log y_{1} + a_{2}(y_{1}- y_{2}) + 1\ ,\ \xi_{2}= \log y_{2}
+ a_{1} (y_{2}- y_{1}) + 1. $$
The defining condition for $u$ is
$du=\xi_{1} dx_{1} + \xi_{2} dx_{1}$ which, after some cancellation
using the fact that \begin{equation}dx_{i} = dy_{i}+ a_{i}( y_{1} dy_{2} + y_{2} dy_{1})\end{equation}  is the differential
$$ \left(\log y_{1} dx_{1} + \log y_{2} dx_{2}\right) + dx_{1} + dx_{2} + (y_{1}-
y_{2})(a_{2} dy_{1}- a_{1} dy_{2}) . $$
Set $V= x_{1} \log y_{1} + x_{2} \log y_{2}$ so
$$ dV= \log y_{1} dx_{1} + \log y_{2} dx_{2} + \frac{x_{1}}{y_{1}} dy_{1}
+ \frac{x_{2}}{y_{2}} dy_{2}. $$
Then $$ du-dV= dx_{1} + dx_{2} + (1-a_{1} y_{2}) dy_{1} + (1-a_{2} y_{2})
dy_{2} + (y_{1}-y_{2}) (a_{2} dy_{1} - a_{1} dy_{2}). $$
Using (62), expanding and cancelling terms,
 one finds that
$$  du-dV= a_{2} y_{1} dy_{1}+ a_{1} y_{2} dy_{2}, $$
so we can take
$$  u= V+ \frac{1}{2} (a_{2} y_{1}^{2} + a_{1} y_{2}^{2}). $$

\

 So far, we have been working rather formally---ignoring the precise domains
 of our functions---and now we return to a global point of view.
 One can verify that, for any $a_{i}> 0$, the map $$(y_{1}, y_{2}) \mapsto (y_{1}+ a_{1} y_{1} y_{2}, y_{2} + a_{2} y_{1} y_{2})$$ yields a diffeomorphism
from the quadrant $\{y_{i}\geq 0\}$ (regarded as a manifold with a corner)
to itself. So we have an inverse diffeomorphism given by functions
$y_{i}=y_{i}(x_{1}, x_{2})$, which can be given by explicit formulae, as
below.
Now define
\begin{equation}u(x_{1}, x_{2})= x_{1} \log y_{1} + x_{2} \log y_{2} + \frac{1}{2}( a_{2}
y_{1}^{2} + a_{1} y_{2}^{2}). \end{equation}

Then $u$ is a convex function on the quarter-plane, satisfying Guillemin boundary
conditions along the axes and $u_{ij}^{ij}=0$. 

\

\

Multiplying $a_{1},a_{2}$ by the same non-zero factor does not change the metric up
to isometry. When $a_{1}=a_{2}$  we have $\xi_{1}+\xi_{2}= \log r +2$ and
the metric is Ricci flat, given by the Gibbons-Hawking construction using
the harmonic function $\frac{1}{\vert r\vert}+1$ on $\bR^{3}$. (We need to
make a linear change of co-ordinates to $\xi_{s}, \xi_{t}$ to fit in with
the discussion above.) This is a standard description of the recognize the {\it Taub-NUT}
metric on $\bR^{4}$, which is well-known to be complete, with curvature in
$L^{2}$. When $a_{1}\neq a_{2}$ the metric is not Ricci-flat. It is easy
to see that it is complete: the author expects, but has not yet checked in
detail, that the curvature is in $L^{2}$. (The definitions above make sense when one of the
$a_{i}$ is zero, and we still get a metric on $\bR^{4}$. But in this case
the curvature is definitely {\it not} in $L^{2}$, so we exclude it.)

We want to discuss the asymptotic behaviour of one of these solutions for
large
$\ux=(x_{1}, x_{2})$. It is convenient to make a linear change of variable
$$   \sigma= a_{2} x_{1}- a_{1} x_{2}\ \ \ , \ \ \tau=a_{1} x_{2}+a_{2} x_{1}.
$$
(So in the Taub-NUT case, when $a_{1}=a_{2}=1/2$ these coincide with the co-ordinates
$s,t$ we used in Section 5.) Then we can solve for $y_{i}$ to find
\begin{equation}  2 y_{1}=(\sigma-1) + \sqrt{ \sigma^{2} + \frac{\tau}{a_{2}}+ 1}, 2y_{2}=
(1-\sigma)+ \sqrt{ \sigma^{2}+ \frac{\tau}{a_{1}}+ 1}. \end{equation}
We consider the behaviour when $\tau$ is large in the three sectors
$\sigma>\epsilon \tau, \vert \sigma \vert \leq \tau, \sigma<-\epsilon \tau$,
for fixed  $\epsilon$.
 If  $\sigma>\epsilon \tau$ we have $y_{1}\sim \sigma$
and $y_{2}=O(1)$ while if $\sigma<-\epsilon \tau$ we have $y_{1}=O(1)$ and
$y_{2}\sim \sigma$. The asymptotics in a sector $-\epsilon \tau\leq \sigma<\epsilon
\tau$ are more complicated, but on the line $\sigma=0$ we have $y_{i}=O(\sqrt{\tau})$.
Now substituting in the formula (63) for $u$ we find that when $\sigma=0$, $u=
O(\tau\log \tau)$. This is the same growth rate as in the Euclidean case.
When $\sigma>\epsilon \tau$ we have $u\sim \frac{a_{1}}{2} \sigma^{2}$ while
if $\sigma<-\epsilon \tau$ we have $u\sim \frac{a_{2}}{2} \sigma^{2}$. Thus
$u$ grows much faster away from the line $\sigma=0$ than it does along this line,
but the growth rates in the two regions $\pm \sigma \geq \epsilon \tau$ are different.
\subsection{Discussion}

\

Suppose we have a convergent sequence of data sets $(P^{(\alpha)},
\sigma^{(\alpha)}, A^{(\alpha)})$, with solutions $u^{(\alpha)}$ but the limit
$(P^{(\infty)}, \sigma^{(\infty)}, A^{(\infty)})$ does {\it not} satisfy the
positivity condition on $L=L_{A^{(\infty)}, \sigma^{(\infty)}}$ discussed in the Introduction. How do the solutions $u^{(\alpha)}$ behave as $\alpha\rightarrow
\infty$? By the results of \cite{kn:D1}, there is an affine-linear function $\lambda$
which changes sign on $P$, such that $L(\lambda^{+})=0$, where  $\lambda^{+}=\max(\lambda,0)$.
Suppose for the moment that this is the unique function with this property,
up to a factor. 
 The line (or \lq\lq crease'') $\{\lambda=0\}$ divides $P$ into two pieces.
   What we expect is that on the interior of each piece the $u^{(\alpha)}$
 converge, after suitable normalisation,  but the normalisations required
 are different, and if we normalise $u^{\alpha}$ on the region $\{\lambda<0\}$
 then on the other region $\{\lambda>0\}$ the functions blow up as
 $$   u^{\alpha} \sim n_{\alpha} \lambda^{+} $$
  with scalars $n_{\alpha}\rightarrow \infty$.
 
  We can easily
 write
 down explicit examples of this behaviour, which is essentially a one-dimensional
 phenomenon. Fix a family of even functions $f_{\epsilon}$  on the interval
 $[-1,1]$, parametrised by $\epsilon\in [0,1]$, with $f_{\epsilon}(x)= x^{2}+\epsilon^{2}$ for $\vert x\vert\leq
 1/2$, with  $f_{\epsilon}(x)= 1-\vert x \vert$ for $\vert x\vert$ close to
 $1$ and with $f_{\epsilon}(x)>0$ except when $\vert x\vert=1$ or $x=\epsilon=0$. We can obviously do this in such a way the family is smooth in both
 variables. Set $a_{\epsilon}= f_{\epsilon}''$. Then for $\epsilon>0$ the
 one-dimensional version of (1), which is
 \begin{equation}  \left( \frac{1}{ U(x)''}\right)= - a_{\epsilon}, \end{equation}
 has a  solution $U_{\epsilon}$ which is smooth in $(-1,1)$ and satisfies
 the Guillemin boundary condition at $\pm 1$. We just take $U_{\epsilon}$
 to be the solution of the elementary equation
 $$    U_{\epsilon}''= f_{\epsilon}^{-1}, $$
 normalised to $U_{\epsilon}(0)=U'_{\epsilon}(0)=0$, say. Equally obviously,
 the family $U_{\epsilon}$ is unbounded on any neighbourhood of $0$, as $\epsilon\rightarrow
 0$, because the derivative $U'_{\epsilon}$ has limit $\vert x\vert ^{-1}$ which is not
 integrable. If, on the other hand, we define $U_{+,\epsilon}$ to be the
 solution of (65) with $U_{+,\epsilon}(1/2)=U_{+,\epsilon}'(1/2)=0$, then the
 $U_{+,\epsilon}$ converge as $\epsilon\rightarrow 0$ on the interval $(0,1]$.
 Similarly there is a  family $U_{-,\epsilon}$ normalised at $-1/2$ and converging
 over $[-1,0)$. For $\epsilon>0$ we have
 $$   U_{+,\epsilon}(x)-U_{-,\epsilon}(x)= n_{\epsilon} x, $$
 where $n_{\epsilon}\rightarrow \infty$ (and in fact $n_{\epsilon}\sim \log \epsilon^{-1}$). 
To get a two-dimensional example we can simply take $P$ to be the square
$[-1,1]^{2}\subset \bR^{2}$ and $A_{\epsilon}(x_{1}, x_{2})= a_{\epsilon}+1$.
Then for $\epsilon>0$ there is a solution $u_{\epsilon}(x_{1}, x_{2})= U_{\epsilon}(x_{1})+
V(x_{2})$, where $V$ is the symplectic potential for the round metric on
the 2-sphere. In terms of Riemannian geometry  in four-dimensions
described by this family: in the region corresponding to $[-1/2,1/2]\times
[-1,1]$ the $4$-manifolds have the form
$$  H_{\epsilon}\times S^{2}$$
with the product of the round metric of curvature $+1$ on $S^{2}$ and a metric of curvature $-1$ on $H_{\epsilon}$,  approaching a pair of \lq\lq cusps'' as
$\epsilon\rightarrow 0$. The curvature of these metrics is bounded,
uniformly in $\epsilon$, but the diameter tends to infinity and the injectivity
radius to $0$.

This gives a model, albeit somewhat conjectural, for the behaviour  near a \lq\lq crease'' $\{\lambda=0\}$ in the general case, {\it provided that
this crease does not pass through a vertex of $P$}. Suppose on the other
hand that we are in this situation, and take the standard model with the vertex the origin and $P$  coinciding locally with the quarter-plane $\{x_{i}>0\}$.
The crease is a line $ x_{2}= b x_{1}$, where $b>0$. What seems 
likely to be true is that, near the origin, the solutions $u^{(\alpha)}$ are
modelled on the zero scalar curvature metric discussed above with parameters
$a_{1}=1, a_{2}= b $, then scaled down by a factors $r_{\alpha}$, with $r_{\alpha}\rightarrow
0$ as $\alpha\rightarrow \infty$ (which is the same as taking parameters
$a_{1}= r_{\alpha}, a_{2}= b r_{\alpha}$. In other words, we expect these
model solutions to appear as the \lq\lq blow-up'' limits. This picture is
consistent with the  discussion in the preceding subsection of the asymptotics
of the model solutions. Notice that (if this picture is correct) then in
this situation the curvature of the $u^{(\alpha)}$  is not bounded uniformly
in the
family, in contrast to the previous case.

Of course we can envisage somewhat more complicated situations in which we have several
\lq\lq creases'', dividing the polygon into  more parts. This is
discussed in \cite{kn:D1}, and taken much further by Szelyhidi in \cite{kn:Sz}.  Szekelyhidi's work also suggests very strongly that a similar picture  holds for the limiting
behaviour of the Calabi flow. Notice also that this (conjectural) picture
is quite in line with the more general situation, in four-dimensional Riemannian
geometry, described by Anderson in
\cite{kn:An}. 

These models are also useful in understanding the issues involved in the
existence proofs.  Consider the first \lq \lq product'' example, but apply
an affine transformation so the domain $P$ is now a parallelogram, say
$2>x_{1}>0, \vert x_{1}- x_{2}\vert \leq 1$, and the crease is the line
$x_{1}=x_{2}$. Then it is clear that, in the family parametrised by $\epsilon$,
the quantities $D(p)$ are not uniformly bounded, for $p$ in an arbitrarily
small neighbourhood of the origin. As we explained in Section 3, the essential difficulty in the proof of Theorem 3 is to show
that we cannot have a family behaving in this way close to $(0,0)$ unless
the $D(p)$ are also large for some large $p$. The same discussion  applies
at the vertices. It is clear from our description of the asymptotics of the
model solutions that for them the quantity $E(t)$ is unbounded as $t\rightarrow
\infty$. Thus, in the setting of Section 5, we cannot obtain an {\it a priori} bound
on $\Emax$ by \lq\lq local'' considerations around the vertex. The essential
difficulty in the proof in Section 5 is to show that we cannot have a situation where
$E(t)$ is large for some range of $t$ but nevertheless $E(t)$ is bounded
for very large $t$: in particular we cannot have a blow-up of the curvature,
modelled on our explicit solutions, unless the $u^{(\alpha)}$ are already
unbounded in the interior of $P$.

\end{document}

%% file: extremal2a.bbl
\begin{thebibliography}{99}


\bibitem{kn:An}  Anderson, M.T. {\em Canonical metrics on 3-manifolds and
4-manifolds} Asian Jour. Math. {\bf 10} 2006 127-163
\bibitem{kn:CP} Calderbank, D.M.J. and Pedersen, H. {\em Self-dual Einstein
metrics with torus symmetry} Jour. Diff. Geom. {\bf 60} 485-521 2002
\bibitem{kn:D1}  Donaldson, S.K.  {\em Scalar curvature and stability of toric
varieties} Jour. Differential Geometry {\bf 62} 289-349  2002
\bibitem{kn:D2}  Donaldson, S.K.  {\em Interior estimates for solutions of
Abreu's equation} Collectanea Math. {\bf 56} 103-142  2005
\bibitem{kn:D3}  Donaldson, S.K.  {\em Extremal metrics on toric surfaces:
a continuity method} To appear in {\it Jour. Differential Geometry} 
\bibitem{kn:D4} Donaldson, S.K. {\em A generalised Joyce construction for
a family of nonlinear partial differential equations} Preprint
\bibitem{kn:J} Joyce,D.D. {\em Explicit construction of self-dual 4-manifolds}
Duke Math. J. {\bf 77} 519-552 1995
\bibitem{kn:Sz} Szekelyhidi, G. {\em Optimal test configurations for toric
varieties} arxiv 07092687
\bibitem{kn:ZZ}  Zhou, B. and  Zhu, X. {\em A note on K-stability of toric manifolds} arxiv 07060505
\end{thebibliography}
